\documentclass{article}
\usepackage{amssymb}
\usepackage{amsmath}
\usepackage{graphicx}
\usepackage[ruled,vlined]{algorithm2e}
\usepackage{supertabular}
\newtheorem{theorem}{Theorem}
\newtheorem{definition}{Definition}
\newtheorem{lemma}[theorem]{Lemma}
\newtheorem{rem}[theorem]{Remark}

\def\oarrow{{\mathrel{\longrightarrow\mkern-25mu\circ}\;\;}}

\newcommand{\inte }{{\mathrm{int}}\,}

\newcommand{\dom }{\,{\mathrm{dom}}\,}

\newcommand{\Id }{\mbox{Id} \,}
\newcommand{\Fix }{\mbox{Fix} \,}
\newcommand{\Inv }{\mbox{Inv} \,}

\newcommand{\comment}[1]{\mbox{}}

\def\qed{{\hfill{\vrule height5pt width3pt depth0pt}\medskip}}

\begin{document}
\begin{center} {\bf \LARGE  Period doubling  in the R\"ossler system - a computer assisted proof} \\
\vskip 0.5cm
  \vskip\baselineskip
    {\large
    Daniel Wilczak\footnote{
            Research supported by an annual national scholarship for young scientists from
            the Foundation for Polish Science
        },
    Piotr Zgliczy\'nski\footnote{
            Research supported in part by Polish State Ministry of
            Science and Information Technology  grant N201 024 31/2163 
         }
    } \\
 Jagiellonian University, Institute of Computer Science, \\
 Nawojki 11, 30--072  Krak\'ow, Poland \\ e-mail: wilczak@ii.uj.edu.pl, umzglicz@cyf-kr.edu.pl
\vskip\baselineskip \today

\end{center}

\section{Introduction}
The goal of this paper is to show how to produce a piece of
rigorous bifurcation diagram of periodic orbits  for an ODE. We
study the R\"ossler system \cite{R}, one of the  textbook examples
of ODEs generating nontrivial dynamics, for the parameter range
containing two period doubling bifurcations.

According to the discussion in Kuzniecov textbook \cite[Section
2.7]{Ku} there are two extremes in studying bifurcations in
dynamical systems. The first one, going back to Poincar\'e, is to
analyze the appearance (branching) of new invariant objects
(equilibria or periodic orbits) from the known ones as parameters
of the system vary. A good reference for this approach is a
textbook by Chow and Hale~\cite{CH}. On the other extreme, it is
the approach going back to Andronov \cite{An} and Thom \cite{T},
is to study rearrangements (bifurcations)  of the whole phase
portrait under variations of parameters. It is apparent that the
first approach is necessarily one of the initial steps in
attempting to describe the bifurcations in Andronov-Thom sense. In
fact in many dimensional systems (even for planar maps like the
H\'enon map) achieving the complete description of the phase space
portrait and its changes appears to be hopeless  in view of the
results on  the Henon-like maps
\cite{MoraViana,BenedicksCarleson,WangYoung1,WangYoung2}.

While there exists a vast literature on the bifurcation theory,
see for example \cite{AAIS,CH,G,Ku} and references given there,
and also a lot of numerical bifurcation diagrams for various
systems can be found in literature (see for example references in
\cite{Ku}), there are virtually no rigorous results on
bifurcations of periodic orbits for ODEs in dimension three or
higher in the situation, when the periodic orbit undergoing the
bifurcation is not given to us analytically due to some special
symmetries of the system. The basic reason for this is: while
numerical experiments and/or  normal form  computations may
clearly show what is happening (in terms of the bifurcations) we
usually lack any reasonable rigorous estimates about the observed
orbits, which prevents us to turn these observations into rigorous
statements. To obtain the necessary estimates one  needs to
integrate the variational equations describing the partial
derivatives with respect to initial conditions up to order 3 or
higher. This is usually a serious problem for rigorous ODE
solvers. It turns out that the naive approach: applying an ODE
solver to the system of variational equations  does not work,
because the methods dealing with the wrapping effect  used in the
Lohner-type algorithms (the most effective rigorous ODE solvers)
\cite{Lo,Z1,NJ} break down for such system. As the solution of
this problem $C^r$-Lohner algorithm has been proposed in
\cite{WZn} and it is used in the present work.

Concerning the content of the paper regarding   the bifurcation
theory itself, we were forced to reformulate some well known
theorems to make them amenable to computer assisted proofs. It is
a common feature of all bifurcations theorems that the bifurcation
point (or rather a candidate) and all necessary data like the
spectrum and maybe some higher order terms  are always given as
part of the assumptions. But in a nonlinear system we usually do
not have explicitly these data, in fact the existence of the
bifurcation point has to be proved by looking on the behavior of
the system in some neighborhood. This forces us to reformulate
some bifurcation theorems in a semi-local way, we have to
investigate properties of solutions of implicit equations, which
are degenerate (due to the presence of bifurcations). This is the
reason, why from various approaches to bifurcations we chose the
one developed in \cite{CH} and which is based on the
Liapunov-Schmidt reduction.

In our work we focus on the period doubling bifurcation of
periodic orbits for R\"ossler equations, in fact we study the
Poincar\'e map for R\"ossler system. The paper is organized  as
follows: in Sections \ref{sec:impfcond}, \ref{sec:contfp} and
\ref{sec:dyninfo} we discuss the main tools used to produce a
validated piece of the bifurcation diagram containing the period
doubling bifurcations. In the remaining sections we give some
details concerning our results for R\"ossler system.

\section{Basic definitions}

By $\mathbb{N}$, $\mathbb{Z}$, $\mathbb{Q}$, $\mathbb{R}$,
$\mathbb{C}$ we denote the set of natural, integer, rational, real
and complex numbers, respectively. $\mathbb{Z}_-$ and
$\mathbb{Z}_+$ are negative and positive integers,  respectively.
By $S^1$ we will  denote a unit circle on the complex plane.

For $\mathbb{R}^n$  we will denote the norm of $x$ by $\|x\|$ and
if the formula for the norm is not specified in some context, then
it means that one  can use any norm  there. Let $x_0 \in
\mathbb{R}^s$, then $B_s(x_0,r)=\{z \in \mathbb{R}^s \: | \: \|x_0
- z \| < r \}$ and $B_s=B_s(0,1)$.

Let $A:\mathbb{R}^n \to \mathbb{R}^n$ be a linear map. By
$\text{Sp}(A)$ we denote the spectrum of $A$, which is the set of
$\lambda \in \mathbb{C}$, such that there exists $x\neq
\mathbb{C}^n \setminus \{0\}$, such that $Ax=\lambda x$.

For a map $f:X \to Y$ by $\dom(f)$ we will denote the domain of
$f$. For a map $F:X \to X$ we will denote the fixed point set by
$\Fix(F,U)=\{x\in U\ | \ F(x)=x\}$.

Let $x=(x_1,\dots,x_n) \in \mathbb{R}^n$. By $\pi_{i}$ we will
denote the projection on $i$-th coordinate, i.e. $\pi_i(x)=x_i$.
Analogously for any multiindex $\alpha=(i_1,i_2,\dots,i_k)\in\mathbb
Z_+^k$ we define $\pi_\alpha(x)=(x_{i_1},x_{i_2},\dots,x_{i_k})$.
Sometimes the points in the phase space will have coordinates
denoted by different letters, for example $z=(\nu,x,y)$, then we
will index the projection by the names of variables, i.e.
$\pi_{(\nu,x)}(z)=(\nu,x)$ etc.

\begin{definition}
\label{def:fphypmap} Let $f:\mathbb{R}^n \supset \dom(f)\to
\mathbb{R}^n$ be $C^1$. Let $z_0 \in \dom(f)$. We say that
$z_0$ is a \emph{hyperbolic fixed point for $f$} iff $f(z_0)=z_0$
and $\mbox{Sp}(Df(z_0)) \cap S^1=\emptyset$, where $Df(z_0)$ is
the derivative of $f$ at $z_0$.
\end{definition}

\begin{definition}
Consider a map $f:X \supset \dom(f) \to X$. Let $x \in X$. Any
sequence $\{x_k\}_{k \in I}$, where $I \subset \mathbb{Z}$ is a set
containing $0$ and for any  $l_1 < l_2 < l_3$ in $\mathbb{Z}$ if
$l_1,l_3 \in I$, then $l_2 \in I$, such that
\begin{equation*}
  x_0=x, \qquad f(x_i) =x_{i+1}, \qquad \mbox{for $i,i+1 \in I$}
\end{equation*}
will be called  an orbit through $x$. If $I=\mathbb{Z}_-\cup\{0\}$,
then we will say that $\{x_k\}_{k \in I}$ is a full backward orbit
through $x$.
\end{definition}

\begin{definition}
Let $X$ be a topological space and  let the map $f:X \supset
\dom(f) \to X$ be continuous.

 Let $Z \subset \mathbb{R}^n$, $x_0 \in Z$, $Z \subset
\dom(f)$. We define
\begin{eqnarray*}
  W^s_Z(z_0,f)&=&\{ z  \: | \: \forall_{n \geq 0} f^n(z)\in Z, \quad \lim_{n\to\infty} f^n(z)=z_0   \} \\
 W^u_Z(z_0,f)&=&\{ z  \: | \:  \exists \mbox{ $\{x_n\} \subset Z$  a full backward orbit through $z$, such that } \\
     & & \quad \lim_{n\to-\infty} x_n=z_0   \} \\
  W^s(z_0,f)&=&\{ z  \: | \:  \lim_{n \to \infty} f^n(z)=z_0   \} \\
 W^u(z_0,f)&=& \{ z  \: | \:  \exists \mbox{ $\{x_n\}$  a full backward orbit through $z$, such that } \\
     & & \quad \lim_{n\to-\infty} x_n=z_0   \} \\
  \Inv^+(Z,f)&=&\{ z  \: | \: \forall_{n \geq 0} f^n(z)\in Z\} \\
  \Inv^-(Z,f)&=&\left\{ z  \: | \: \exists \mbox{ $\{x_n\} \subset Z$  a full backward orbit through $z$ }
  \right\} \\
  \Inv(Z,f)&=&  \Inv^+(Z,f) \cap  \Inv^-(Z,f)
\end{eqnarray*}
If $f$ is known from the context, then we will usually drop it and
use $W^s(z_0)$, $W_Z^s(z_0)$ etc instead.
\end{definition}

\begin{definition}
\label{def:perdoubl} Let $P_\nu: \mathbb{R}^n \to \mathbb{R}^n$,
where $\nu$ belongs to some interval. We say that $P_\nu$ has
\emph{a period doubling bifurcation} at $(\nu_0,x_0)$ iff there
exists $V=[\nu_1,\nu_2] \times X \subset \mathbb{R} \times
\mathbb{R}^n$, such that the following conditions are satisfied
\begin{itemize}
\item $(\nu_0,x_0) \in \inte V$, $P_{\nu_0}(x_0)=x_0$
\item there exists a continuous function $x_{fp}:[\nu_1,\nu_2]\to \inte X$,
such that
\begin{equation*}
 \Fix(P_\nu,X)=\{ x_{fp}(\nu) \}
\end{equation*}
\item  there exist two continuous curves
$c_i:[\nu_0,\nu_2] \to \inte X$, $i=1,2$ , such that for $\nu \in
[\nu_0,\nu_2]$ holds
\begin{eqnarray*}
c_1(\nu_0)&=&c_2(\nu_0)=x_{fp}(\nu_0) \\
   c_1(\nu) &\neq& c_2(\nu), \quad \nu \neq \nu_0 \\
   P_\nu(c_1(\nu))&=&c_2(\nu), \quad  P_\nu(c_2(\nu))=c_1(\nu)
    \\
   Fix(P^2_\nu,X)&=& \{c_1(\nu), c_2(\nu), x_{fp}(\nu)\}
\end{eqnarray*}
\item the dynamics:

for $\nu \leq  \nu_0$
\begin{displaymath}
  \Inv(X,P_\nu)=\{ x_{fp}(\nu) \}
\end{displaymath}

For $\nu > \nu_0$ the maximal invariant set in $X$ $\Inv(X,P_\nu)$
is equal to
\begin{equation*}
\overline{W_X^u}(x_{fp}(\nu),P) \cap
\left(\overline{W_X^s}(c_1(\nu),P_\nu^2) \cup
\overline{W_X^s}(c_2(\nu) ,P_\nu^2) \right)
\end{equation*}
and is a one-dimensional connected manifold with boundary points $c_1(\nu)$, $c_2(\nu)$.
\end{itemize}
\end{definition}

\section{Derivation of the conditions for the occurrence of
 the period doubling bifurcation}
\label{sec:impfcond} The goal of this section is to present the
set of conditions, which guarantee the existence of period
doubling bifurcation for a given map, and which can be verified
using rigorous numerics. The main tools used are the
Liapunov-Schmidt reduction \cite{CH} and the implicit function
theorem.

Assume that we have a parameter dependent map $z \mapsto
P(\nu,z)$, which apparently undergoes the period doubling
bifurcation as the parameter $\nu$ changes. Let $z_{fp}(\nu)$ be a
fixed point curve for $P$. We assume that it is regular and we can
compute it and its all derivatives.

To prove the existence of the period doubling bifurcation we
proceed as in \cite{CH}. First we perform the Liapunov-Schmidt
reduction to obtain a function $G:\mathbb{R} \times \mathbb{R}
\supset \dom(G) \to \mathbb{R}$, whose zeros correspond to fixed
points and period two points of $P_\nu$ and then we try to
describe the solution set for equation $G(\nu,x)=0$. Next, through
some additional computation of eigenvalues we will be able to
decide about the hyperbolicity of bifurcating periodic orbits.

The basic steps of the Liapunov-Schmidt reduction for $P^2$ are:
\begin{itemize}
  \item  we choose good coordinates $(x,y) \in \mathbb{R} \times
  \mathbb{R}^{n-1}$.
  It is desirable to choose $x$ in the
  approximate bifurcation direction (in the eigendirection corresponding to $-1$
   eigenvalue at the bifurcation point).
 \item let $Z=[\nu_1,\nu_2] \times [x_1,x_2]$ and  $Y \subset \mathbb{R}^{n-1}$
  be such that the apparent bifurcation point $(\nu_0,x_0,y_0)$
  belongs to the interior of $Z \times Y$
  \item we need to show that there exists a  function $y(\nu,x)$, defined on $Z$
    with the values in $Y$,  such that
  \begin{equation}
     y - \pi_y (P^2_\nu(x,y)) =0 \quad \mbox{for $(\nu,x,y) \in Z \times Y$  iff}
    \quad y=y(\nu,x).  \label{eq:for-y}
  \end{equation}

 \item the bifurcation function $G:Z \to \mathbb R$ is defined by
 \begin{equation}
   G(\nu,x)=x - \pi_x (P^2_\nu(x,y(\nu,x))). \label{eq:bif1}
 \end{equation}
\end{itemize}

Now, we have to find the solution set of the following equation
\begin{equation}
  G(\nu,x)=0, \qquad (\nu,x)\in Z. \label{eq:bif2}
\end{equation}

Let $x_{fp}(\nu)=\pi_x(z_{fp}(\nu))$ be the $x$-coordinate of the
fixed point curve. We assume that $[\nu_1,\nu_2] \subset
\dom(x_{fp})$ and $x_{fp}([\nu_1,\nu_2]) \subset [x_1,x_2]$.
Therefore we have
\begin{equation}
  G(\nu,x_{fp}(\nu))=0. \label{eq:fixp}
\end{equation}

The idea of solving (\ref{eq:bif2}) goes as follows: we introduce
a new bifurcation function
\begin{equation}
g(\nu,x)=\frac{G(\nu,x)}{x-x_{fp}(\nu)} \label{eq:gmGD}
\end{equation}
 and then we solve equation $g(\nu,x)=0$ by the implicit function theorem.

Observe that expression (\ref{eq:gmGD}) defining  $g(\nu,x)$
contains zero in the denominator, moreover usually the exact value
of $x_{fp}(\nu)$ is not known, therefore the formula
(\ref{eq:gmGD}) appears to be useless in rigorous computations.
The next lemma will give  us an integral representation of $g$,
which will not contain any singularities and therefore it is well
suited for rigorous numerics.

\begin{lemma}
Assume $F:\mathbb{R}^n \to \mathbb{R}^s$ is $C^1$.  Let $x,y \in
\mathbb{R}^n$. Then
\begin{equation*}
  F(x) - F(y)=\int_0^1 \frac{\partial F}{\partial x}(t(x-y) +
  y)dt \cdot (x - y)
\end{equation*}
\end{lemma}

Hence we can define equivalently $g:[\nu_1,\nu_2]\to [x_1,x_2]$ by
\begin{equation}
  g(\nu,x)= \int_0^1 \frac{\partial G}{\partial x}(\nu,t(x-x_{fp}(\nu)) +
  x_{fp}(\nu))dt. \label{eq:def-g}
\end{equation}
We obtain
\begin{equation*}
  G(\nu,x)=(x-x_{fp}(\nu))g(\nu,x).
\end{equation*}

Therefore, we have to determine the solution set of  the following
equation
\begin{equation}
  g(\nu,x)=0 \quad (\nu,x) \in Z,  \label{eq:g-bif}
\end{equation}
where $g$ is defined in (\ref{eq:def-g}).

 In the case of the period doubling bifurcation we expect
solutions of (\ref{eq:g-bif}) to form a regular curve. The
following lemma gives a set of conditions, which implies this
fact.

\begin{lemma}
\label{lem:saddle-node} Let $Z=[\nu_1,\nu_2] \times [x_1,x_2]$.
Assume that $g:Z \to \mathbb{R}$ is a $C^k$-function, $k \geq 2$.

Assume that
\begin{eqnarray}
   \frac{\partial^2 g}{\partial x^2}(Z) &>& 0, \label{eq:d2g} \\
    \frac{\partial g}{\partial \nu}(Z) &<& 0, \label{eq:dg}  \\
   g(\nu_1,x) &>&0, \qquad \mbox{for $x \in [x_1,x_2]$} \label{eq:gbd1} \\
   g(\nu_2,x_1) &>& 0,  \label{eq:gbd2} \\
    g(\nu_2,x_2) &>& 0, \label{eq:gbd2.2} \\
   g(\nu_2,x_0) &<& 0, \qquad \mbox{for some $x_0 \in (x_1,x_2)$} \label{eq:gbd3}
\end{eqnarray}
Then there exist $\bar{x}_1,\bar{x}_2$, such that $x_1 < \bar{x}_1
< x_0 < \bar{x}_2 < x_2$ and there exists a function
$\nu:[\bar{x}_1,\bar{x}_2] \to [\nu_1,\nu_2]$ of class $C^k$, such
that
\begin{equation*}
  \{(\nu,x) \in Z \ | \  g(\nu,x)=0 \ \}=\{(\nu(x),x), x\in
  [\bar{x}_1,\bar{x}_2]\}.
\end{equation*}

 Moreover, there exists $\bar{x}_0 \in
(\bar{x}_1,\bar{x}_2)$ such that
\begin{eqnarray*}
  \nu'(x) &>& 0, \qquad x \in (\bar{x}_0,\bar{x}_2) \\
  \nu'(x) &<& 0, \qquad x \in (\bar{x}_1, \bar{x}_0) \\
  \nu(x) &>& \nu_1, \qquad x \in [\bar{x}_1,\bar{x}_2] \\
  \nu(\bar{x}_1) &=&  \nu(\bar{x}_2)=\nu_2.
\end{eqnarray*}
\end{lemma}
\textbf{Proof:} Observe first that from condition (\ref{eq:d2g})
it follows that for any given $\nu \in [\nu_1,\nu_2]$ and any $c
\in \mathbb{R}$ the equation
\begin{equation*}
  g(\nu,x)=c,
\end{equation*}
has at most two solutions in $[x_1,x_2]$.

From this observation and equations (\ref{eq:gbd2}--\ref{eq:gbd3})
it follows that there exist  $\bar{x}_1$ and $\bar{x}_2$, such
that
\begin{eqnarray*}
  x_1\ <\ \bar{x}_1\ <& x_0& <\ \bar{x}_2\ <\ x_2 \\
  \{x \in [x_1,x_2]\ | \ g(\nu_2,x)=0 \ \} &=& \{\bar x_1,\bar x_2\} \\
  g(\nu_2,x) &>& 0, \quad \mbox{for $x < \bar{x}_1$ or $x >
  \bar{x}_2$}\\
   g(\nu_2,x) &<& 0, \quad \mbox{for $x \in
   (\bar{x}_1,\bar{x}_2)$}.
\end{eqnarray*}
From the above conditions and conditions (\ref{eq:dg}) and
(\ref{eq:gbd1}) it follows immediately, that there exists function
$\nu:[\bar{x}_1,\bar{x}_2] \to [\nu_1,\nu_2]$, such that
\begin{equation*}
  \{(\nu,x) \in Z \ | \  g(\nu,x)=0 \ \}=\{(\nu(x),x), x\in
  [\bar{x}_1,\bar{x}_2]\}.
\end{equation*}
By the implicit function theorem   function $\nu(x)$ is of class
$C^k$.

It remains to show the existence of a unique minimum of $\nu(x)$
and its monotonicity properties.

Let $y\in (\bar{x}_1,\bar{x}_2)$ be any critical point of
$\nu(x)$, i.e $\dot{\nu}(y)=0$. We will show that
$\ddot{\nu}(y)>0$.

By differentiating twice equation $g(\nu(x),x)=0$ we obtain
\begin{eqnarray*}
\frac{\partial^2 g}{\partial \nu^2}(\nu(x),x) (\dot{\nu}(x))^2 +
2\frac{\partial^2 g}{\partial \nu \partial x}(\nu(x),x)
\dot{\nu}(x) + \\
\frac{\partial g}{\partial \nu}(\nu(x),x) \ddot{\nu}(x) +
\frac{\partial^2 g}{\partial x^2}(\nu(x),x) =0
\end{eqnarray*}
Therefore for $y$ we have
\begin{eqnarray*}
 0&=&\frac{\partial g}{\partial \nu}(\nu(y),y) \ddot{\nu}(y) +
\frac{\partial^2 g}{\partial x^2}(\nu(y),y)  \\
\ddot{\nu}(y) &=& - \left(\frac{\partial g}{\partial \nu}(\nu(y),y)
\right)^{-1} \frac{\partial^2 g}{\partial x^2}(\nu(y),y)\ >\ 0.
\end{eqnarray*}
We see that all critical points are strong local minima. This
implies that the set of critical points consists from just one
point.
 \qed

The model for Lemma~\ref{lem:saddle-node} is given by the function
$g_1(\nu,x)=x^2 - \nu$ in the neighborhood of point $(0,0)$. By
changing signs of $\nu$ and $g$ we obtain the following model
functions $g_2(\nu,x)=\nu + x^2$, $g_3(\nu,x)=\nu-x^2$ and
$g_4(\nu,x)=-\nu-x^2$ for which we can state analogous lemmas.

Now we can formulate a lemma based on the implicit function
theorem addressing the assumptions implying intersection of curves
solving equation $G(\nu,x)=0$, where $G$ arises in through the
Liapunov-Schmidt reduction in the context of the period doubling
bifurcation.
\begin{lemma}
\label{lem:Gsol}
 Let $Z=[\nu_1,\nu_2] \times
[x_1,x_2]$. Assume that $G:Z \to \mathbb{R}$ is a $C^k$-function,
$k \geq 3$.

Assume that there exists a $C^k$-function $x_{fp}:[\nu_1,\nu_2]
\to (x_1,x_2)$, such that $G(\nu,x_{fp}(\nu))=0$ for $\nu \in
[\nu_1,\nu_2]$.

Assume that
\begin{eqnarray}
 \frac{\partial^3 G}{\partial x^3}(Z) > 0 \label{eq:d3gd3x} \\
 \frac{\partial^2 G}{\partial x \partial \nu} (Z) +
  \frac{\partial^2 G}{\partial x^2} (Z)  x'_{fp}([\nu_1,\nu_2]) \cdot
  [0,1] < 0. \label{eq:d2gdxdnu+}
\end{eqnarray}

We assume that following conditions are satisfied for some $x_1
\leq \delta_1 < x_{fp}(\nu_1) < \delta_2 \leq x_2$
\begin{eqnarray}
   G(\nu_1,[\delta_2,x_2]) > 0, \qquad   G(\nu_1,[x_1,\delta_1])
  &<& 0  \label{eq:sgn+Gnu11} \\
    \frac{\partial G}{\partial x}(\nu_1,[\delta_1,\delta_2]) &>&  0,
  \label{eq:sgn+Gnu12}
  \\
    G(\nu_2,x_2) >0, \qquad   G(\nu_2,x_1) &<& 0  \label{eq:sgn+Gnu21} \\
    \frac{\partial G}{\partial x}(\nu_2,x_{fp}(\nu_2)) &<& 0
  \label{eq:sgn+Gnu22}
\end{eqnarray}

Then there exist $x_1 < \bar{x}_1 < \bar{x}_2 < x_2$, such that
$x_{fp}(\nu_2) \in (\bar{x}_1 , \bar{x}_2)$
 and a
function $\nu:[\bar{x}_1,\bar{x}_2] \to [\nu_1,\nu_2]$ of class
$C^{k-1}$, such that
\begin{multline}
 \left\{(\nu,x) \in Z \ | \  G(\nu,x)=0 \ \right\}= C_{fp}\cup C_{per} = \\
    \left\{ (\nu,x_{fp}(\nu)), \ \nu \in [\nu_1,\nu_2]\right\}
    \cup
    \left\{(\nu(x),x), x\in [\bar{x}_1,\bar{x}_2]\right\}
    \label{eq:Gsol2krzywe}
\end{multline}
and the intersection of curves $C_{fp}$ and $C_{per}$ contains
exactly one point.

 Moreover, there exists $\bar{x}_0 \in
(\bar{x}_1,\bar{x}_2)$ such that
\begin{eqnarray*}
   \nu'(x) &>& 0, \qquad x \in (\bar{x}_0,\bar{x}_2) \\
   \nu'(x) &<& 0, \qquad x \in (\bar{x}_1, \bar{x}_0)
  \\
  \nu(x) &>& \nu_1, \qquad x \in [\bar{x}_1,\bar{x}_2] \\
  \nu(\bar{x}_1) &=&  \nu(\bar{x}_2)=\nu_2.
\end{eqnarray*}

\end{lemma}
\textbf{Proof:} For the proof we want  to apply to
Lemma~\ref{lem:saddle-node}. For this end we define $g$ as in
(\ref{eq:def-g}).

We start by showing that (\ref{eq:d3gd3x}) and
(\ref{eq:d2gdxdnu+}) imply that $ \frac{\partial^2 g}{\partial
x^2}(Z)> 0$ and $ \frac{\partial g}{\partial \nu}(Z)  < 0$,
respectively.

We have
\begin{eqnarray*}
  \frac{\partial^2 g}{\partial x^2}(\nu,x)=\int_0^1 \frac{\partial^3 G}{\partial x^3}(\nu,t(x-x_{fp}(\nu)) +
  x_{fp}(\nu))t^2dt.
\end{eqnarray*}
Hence from (\ref{eq:d3gd3x})  we obtain immediately that $
\frac{\partial^2 g}{\partial x^2}(Z) > 0$.

\begin{multline*}
  \frac{\partial g}{\partial \nu}(\nu,x)=\int_0^1\left( \frac{\partial^2 G}{\partial x \partial \nu}(\nu,t(x-x_{fp}(\nu)) +
  x_{fp}(\nu)) + \right. \\
  \left. \frac{\partial^2 G}{\partial x^2}(\nu,t(x-x_{fp}(\nu)) +
  x_{fp}(\nu))(1-t) x'_{fp}(\nu)  \right)dt \subset \\
 \frac{\partial^2 G}{\partial x \partial \nu} (Z) +
  \frac{\partial^2 G}{\partial x^2} (Z)  x'_{fp}([\nu_1,\nu_2]) \cdot
  [0,1].
\end{multline*}
This and (\ref{eq:d2gdxdnu+}) imply that $ \frac{\partial
g}{\partial \nu}(Z)< 0$.

To obtain condition (\ref{eq:gbd1})  we need to split the interval
$[x_1,x_2]$ into three parts $[x_1,\delta_1]$,
$[\delta_1,\delta_2]$ and $[\delta_2,x_2]$, so that in the middle
part we have the zero of $G(\nu_1,\cdot)$ and we need to use there
the integral representation of $g$. On the remaining parts it is
enough to verify the signs of $G$. Hence we see that conditions
(\ref{eq:sgn+Gnu11}--\ref{eq:sgn+Gnu12}) imply (\ref{eq:gbd1}).

The remaining assumptions in Lemma~\ref{lem:saddle-node} follow
easily from (\ref{eq:sgn+Gnu21}--\ref{eq:sgn+Gnu22}). Now we use
Lemma~\ref{lem:saddle-node} to obtain function $\nu(x)$ and
condition (\ref{eq:Gsol2krzywe}).

It remains to show that  curves $C_{fp}$ and $C_{per}$ defined by
(\ref{eq:Gsol2krzywe}) intersect exactly in one point. Observe
that these curves intersect because curve $C_{fp}$ cuts $Z$ into
two pieces and the end points of the second curve belong to
different components, which follows directly from the fact that
$x_{fp}(\nu_2) \in (\bar{x}_1 , \bar{x}_2)$.

Now we turn to the question of the uniqueness of the intersection
point.

Let $\alpha,\beta \in [\nu_1,\nu_2]$, $\alpha < \beta$. For $t \in
[0,1]$ let $\nu_t=t\alpha + (1-t)\beta$ and
 $x_t=tx_{fp}(\alpha) + (1-t)x_{fp}(\beta)$. Observe that for each
$t \in [0,1]$ point $(\nu_t,x_t)$ belongs to $Z$. Let $\theta \in
(\alpha,\beta)$ be such that
$x'_{fp}(\theta)=\frac{x_{fp}(\alpha)-x_{fp}(\beta)}{\alpha - \beta
}$. We have
\begin{multline*}
  \frac{\partial G}{\partial x}(\alpha,x_{fp}(\alpha)) - \frac{\partial G}{\partial
  x}(\beta,x_{fp}(\beta)) = \\
  \int_0^1 \left( \frac{\partial^2 G}{\partial x \partial
  \nu}(\nu_t,x_t)(\alpha - \beta) +
   \frac{\partial^2 G}{\partial x^2}(\nu_t,x_t)(x_{fp}(\alpha) -  x_{fp}(\beta))\right)dt= \\
 \left( \int_0^1 \left( \frac{\partial^2 G}{\partial x \partial
  \nu}(\nu_t,x_t) +
   \frac{\partial^2 G}{\partial x^2}(\nu_t,x_t)x'_{fp}(\theta)\right)dt \right)
  (\alpha - \beta) \subset \\
   \left(\frac{\partial^2 G}{\partial x \partial \nu} (Z) +
  \frac{\partial^2 G}{\partial x^2} (Z)  x'_{fp}([\nu_1,\nu_2]) \cdot
  [0,1] \right)(\alpha -  \beta).
\end{multline*}
Therefore, from above computations and assumption
(\ref{eq:d2gdxdnu+}) it follows that the function $ \nu \mapsto
\frac{\partial G}{\partial x}(\nu,x_{fp}(\nu)) $ is injective on
$[\nu_1,\nu_2]$. Observe that from (\ref{eq:def-g}) it follows
that, if $(\nu,x_{fp}(\nu))\in C_{fp}\cap C_{per}$ then
$\frac{\partial G}{\partial x}(\nu,x_{fp}(\nu))=0$, so the
intersection of $C_{fp}$ and $C_{per}$ contains at most one point.
\qed

Observe that in the above lemma we cannot make the claim that the
intersection point of the curves, which solve equation $G(\nu,x)=0$
is exactly in $(\nu(\bar{x}_0),\bar{x}_0)$. This can be easily seen
in the following example. Let $G(\nu,x)=(x-1)(x^2 - \nu)$, $x_1=-2$,
$x_2=2$, $\nu_1=-1$ and $\nu_2=1$. It is easy to see that all
assumptions of Lemma~\ref{lem:Gsol} are satisfied, but the
intersection of the curves $(\nu(x)=x^2,x)$ and $(\nu,x(\nu)=1)$ is
not $(0,0)$. On the other hand in the context of the period doubling
bifurcation the intersection point is $(\nu(\bar{x}_0),\bar{x}_0)$,
but we cannot infer such conclusion from Lemma~\ref{lem:Gsol} and we
need to use the information about the dynamical origin of function
$G$. Now we state the theorem which addresses this issue.

\begin{theorem}
\label{thm:per-doubling} Let $P_\nu: \mathbb{R}^n \supset
\dom(P_\nu) \to \mathbb{R}^n $, where $\nu \in I \subset \mathbb{R}$
be one-parameter family of maps of class $C^k$ ($k \geq 3$), both
with respect to the parameter $\nu$ and $x \in \mathbb{R}^n$.

Let $Z=[\nu_1,\nu_2] \times [x_1,x_2]$ and $Y\subset
\mathbb{R}^{n-1}$ be a closure of open set, such that $[x_1,x_2]
\times Y \subset \dom(P^2_\nu) $ for $\nu \in [\nu_1,\nu_2]$.
Assume that
\begin{description}
\item[A1] for any $(\nu,x) \in Z$ there exists a unique $y=y(\nu,x) \in \inte
Y$, such that $y - \pi_y (P^2_\nu(x,y))=0$. Moreover, we assume that
$y:Z \to Y$ is $C^k$.
\item[A2] there exists  $C^k$-function $x_{fp}:[\nu_1,\nu_2]\to (x_1,x_2)$,
such that for $\nu\in[\nu_1,\nu_2]$ holds
\begin{equation}
 \Fix(P_\nu,[x_1,x_2] \times Y)=\{ (x_{fp}(\nu),y(\nu,x_{fp}(\nu)))\}
\end{equation}
\item[A3] Let
\begin{equation*}
  G(\nu,x)=x- \pi_x(P^2_\nu(x,y(\nu,x))), \quad \mbox{for $(\nu,x) \in
  Z$}.
\end{equation*}
Assume that $G$ and $x_{fp}$ satisfy assumptions of
Lemma~\ref{lem:Gsol} and let $\bar{x}_1$,$\bar{x}_2$, $\bar{x}_0$
and $\nu:[\bar{x}_1,\bar{x}_2]\to [\nu_1,\nu_2]$ be as in the
assertion of Lemma~\ref{lem:Gsol}.

\end{description}

Then the fixed point set  of $P^2_\nu$ for $\nu \in [\nu_1,\nu_2]$,
i.e.
\begin{equation*}
\left\{ (\nu,x,y) \in Z \times Y \ | \ P_\nu^2(x,y)=(x,y) \right\}
\end{equation*}
is equal to the sum of the fixed point set for $P_\nu$
\begin{displaymath}
Per_1=\{ (\nu,x_{fp}(\nu),y(\nu,x_{fp}(\nu))) \ | \ \nu \in
[\nu_1,\nu_2] \}
\end{displaymath}
 and the period-$2$ points set
\begin{displaymath}
Per_2=
      \{ (\nu(x),x,y(\nu,x)) \ | \  x\in [\bar{x}_1,\bar{x}_2] \}.
\end{displaymath}

Sets $Per_1$ and $Per_2$ have exactly one common point
$(\nu_b,z_b)$ given by
\begin{equation*}
   (\nu_b,z_b)=(\nu(\bar{x}_0),(\bar{x}_0,y(\nu(\bar{x}_0),\bar{x}_0)).
\end{equation*}

 Moreover, the projections of $Per_1$ and $Per_2$ onto
$(\nu,x)$-plane have exactly one common point $(\nu_b,x_b)$ given
by
\begin{equation*}
   (\nu_b,x_b)=(\nu(\bar{x}_0),\bar{x}_0).
\end{equation*}

\end{theorem}
\textbf{Proof:} From the construction of the bifurcation function
$G$ and our assumptions we immediately obtain that
\begin{equation*}
   \{ (\nu,x,y) \in Z \times Y \ | \
P_\nu^2(x,y)=(x,y) \}=Per_1 \cup Per_2.
\end{equation*}

From Lemma~\ref{lem:Gsol} we know that projections onto
$(\nu,x)$-plane of sets $Per_1$ and $Per_2$ intersect exactly in
one point, say $(\bar{\nu},\bar{x})$. Observe that the point
$(\bar{\nu},\bar{x},y(\bar{\nu},\bar{x}))$ belongs to the
intersection of $Per_1$ and $Per_2$.

It remains to show that
$(\bar{\nu},\bar{x})=(\nu(\bar{x}_0),\bar{x}_0)$. We will show
that the function $x \mapsto \nu(x)$ has a local extremum at
$\bar{x}$. This will imply that $\bar{x}=\bar{x}_0$, because by
Lemma~\ref{lem:Gsol} $\bar{x}_0$ is the only local extremum of
$\nu(x)$.

We reason  by contradiction. Let us assume that $\nu'(\bar{x})
\neq 0 $. Let $U=U_{\nu} \times U_{x} \times U_y$, where $U_\nu
\subset [\nu_1,\nu_2]$, $U_x \subset [x_1,x_2]$ and $U_y \subset
Y$, be neighborhood of $(\bar{\nu},\bar{x},y(\bar{\nu},\bar{x}))$,
such that
\begin{eqnarray}
  P_\nu(x,y) &\in& \inte ([x_1,x_2] \times Y), \quad \mbox{for $(\nu,x,y) \in U$}\nonumber\\
  \nu(a) &\neq& \nu(b), \quad \mbox{for $a,b \in U_x$ and $a\neq
  b$}. \label{eq:1to1}
\end{eqnarray}
Such $U$ exists because $(\bar{x},y(\bar{\nu},\bar{x}))$ is a
fixed point for $P_{\bar{\nu}}$ and $(\bar{x},y(\bar{\nu},\bar{x})) \in
\inte ([x_1,x_2]\times Y)$.

Let us take $v \in U_x$, such that $v \neq \bar{x}$. Then
$(v,y(\nu(v),v))$ is not a fixed point for $P_{\nu(v)}$.  Points
$(v,y(\nu(v),v))$ and   $P_{\nu(v)}(v,y(\nu(v),v))$ are different,
both belong to $Z$ and are period-2 points for $P_{\nu(v)}$.
Therefore they both belong to $Per_2$ and
\begin{equation}
  \nu(\pi_x P_{\nu(v)}(v,y(\nu(v),v)))=\nu(v). \label{eq:nueq}
\end{equation}
Observe that from the continuity it follows that
\begin{equation*}
\lim_{v \to \bar{x}}\pi_x P_{\nu(v)}(v,y(\nu(v),v)) = \pi_x
P_{\nu(\bar{x})}(\bar{x},y(\nu(\bar{x}),\bar{x}))=\bar{x}.
\end{equation*}
From the above observation it follows  that for $v$ sufficiently  close to $\bar{x}$
points $v$ and $\pi_x P_{\nu(v)}(v,y(\nu(v),v))$ are in  $U_x$, but in this situation
condition (\ref{eq:nueq}) contradicts (\ref{eq:1to1}). This proves that
$\bar{x}=\bar{x}_0$.

\qed

\subsection{Hyperbolicity of bifurcating solutions}
\label{subsec:hypbif} The Liapunov-Schmidt projection does not
give any direct information about the dynamical character of the
bifurcating objects. The required information concerning the
hyperbolicity is of course contained in the spectra of $DP_\nu$
and $DP^2_\nu$ and its derivatives.
 Below
we present a lemma addresing this issue.
\begin{lemma}
\label{lem:pd-hyper} Assume that  $P_\nu:\mathbb{R}^n \to
\mathbb{R}^n$ for $\nu \in [\nu_1,\nu_2]$ satisfies all
assumptions of Theorem~\ref{thm:per-doubling} and in the sequel we
will use all the notation introduced there.

Let $z_{fp}(\nu)=(x_{fp}(\nu),y(\nu,x_{fp}(\nu)))$.

\begin{description}
\item[fixed points:]

 Assume that there exists $\epsilon >0$, such that
for all $\nu \in [\nu_1,\nu_2]$ holds
\begin{equation*}
 \mathrm{Sp}\left(\frac{\partial P_\nu}{\partial
 z}\left(z_{fp}(\nu)\right)\right) =   A_\nu \cup B_\nu \cup \{ \lambda(\nu) \},
\end{equation*}
 where $\lambda(\nu) \in \mathbb{R}$   has the multiplicity one,
 $A_\nu \subset \{\alpha \in \mathbb{C},\quad |\alpha| < 1 - \epsilon
\}$
  and $B_\nu \subset  \{\beta \in \mathbb{C},\quad |\beta| > 1 + \epsilon \}$.
  Moreover, we assume
 that
 \begin{eqnarray*}
   \lambda(\nu_1) &\subset& (-1, 1) \\
   \lambda(\nu_2) &<& -1 \\
   \frac{d \lambda}{d \nu}(z_{fp}(\nu)) &<&0, \qquad \nu \in [\nu_1,\nu_2]
 \end{eqnarray*}
Then  the fixed points for $P_\nu$ on curve $z_{fp}(\nu)$ are
hyperbolic for $\nu \in [\nu_1,\nu_2] \setminus
\{\nu(\bar{x}_0)\}$ and
\begin{equation*}
  \dim W^u(z_{fp}(\nu^-),P_{\nu^-}) + 1=
  \dim W^u(z_{fp}(\nu^+),P_{\nu^+})
\end{equation*}
for any $\nu_1 \leq \nu^- < \nu(\bar{x}_0) < \nu^+ \leq \nu_2$.

\item[period-2 points:]

 Assume that there exists $\epsilon >0$, such that on the $Per_2$ curve (i.e. for $x \in [\bar{x}_1,\bar{x}_2]$)
 holds
 \begin{equation*}
\mathrm{Sp}\left(\frac{\partial P^2_{\nu(x)}}{\partial
z}\left(x,y(\nu(x),x)\right)\right)=A_x \cup B_x \cup \{\gamma(x)\}
\end{equation*}
 where $\gamma(x) \in \mathbb{R}$  has the multiplicity one,
 $A_x \subset \{\alpha \in \mathbb{C},\quad |\alpha| < 1 - \epsilon
\}$
  and $B_x \subset  \{\beta \in \mathbb{C},\quad |\beta| > 1 + \epsilon \}$.
  Moreover, we assume
 that
 \begin{eqnarray*}
   \frac{d^2 \gamma}{dx^2}(x) &<& 0, \qquad x \in [\bar{x}_1,\bar{x}_2] \\
    0 < \gamma(\bar{x}_1) &<& 1.
 \end{eqnarray*}
 Then for $x \in [\bar{x}_1,\bar{x}_2] \setminus \{ \bar{x}_0\}$
 the period two points  $z_{d}(x)=(x,y(\nu(x),x))$ for $P_{\nu(x)}$ are
 hyperbolic and
 \begin{eqnarray*}
    \gamma(\bar{x}_0)&=&1 \\
    0 < \gamma(x) &<& 1,  \\
    \dim W^s(z_{d}(x), P^2_{\nu(x)}) &=&
         \dim W^s(z_{fp}(\nu^-),P_{\nu^-})
 \end{eqnarray*}
 for any $\nu_1 \leq \nu^- < \nu(\bar{x}_0)$  and $x \in [\bar{x}_1,\bar{x}_2] \setminus \{\bar{x}_0\}$
\end{description}
\end{lemma}
\textbf{Proof:} The statement about the hyperbolicity of fixed
points is obvious.

 For the proof of the second part it is enough to observe that in
the bifurcation point holds
\begin{equation*}
  \lambda(\nu(\bar{x}_0))=-1, \quad \gamma(x_0)=
  \lambda(\nu(\bar{x}_0))^2=1, \quad \frac{d
  \gamma}{dx}(\bar{x}_0)=0.
\end{equation*}
 \qed

\section{Continuation}
\label{sec:contfp}

To apply the tools described in Section~\ref{sec:impfcond} in the
part regarding the existence of the Liapunov-Schmidt reduction  we
need to prove the existence and uniqueness (locally) of solution
of the equation of the form $f(a,y)=0$ for a given $a$, where $y
\in \mathbb{R}^n$ and $a$ is a parameter. Similarly, when
continuing the fixed point curve or period-2 point curve we have
solve the existence and the local uniqueness of the solution of
$x-P^i(a,x)=0$, where $a$ is the parameter. It turns out that both
of the above mentioned tasks, can be handled by the same tools.

In this section we will discuss such tools, the first one consists
of classical interval analysis tools: interval Newton method
\cite{A,Mo,N} and Krawczyk method \cite{A,Kr,N}, which can be seen
as clever interval versions of the standard Newton method. These
methods work very efficiently in the situation, where the solution
sought is well isolated from other solutions   and it requires
$C^1$-estimates, only. The second approach, which is based on the
implicit function theorem deals with situation, when we are close
to the bifurcation point and therefore there are several solutions
close to one another, as in  the case of the period doubling we
have the fixed point and period two points in a small
neighborhood.

\subsection{Two methods for proving the existence of zeros for a map.}

Let $A\subset \mathbb R^n$. By $[A]_I$ we will denote an
interval enclosure of the set $A$, i.e. the smallest set of the
form $[A]_I=[a_1,b_1]\times \cdots \times[a_n,b_n]$, such that
$A\subset[A]_I$, where $a_i,b_i\in\mathbb R^n\cup\{\pm\infty\}$.

\begin{theorem}(\textbf{Interval Newton Method}
\cite{A,Mo,N})\label{thm:NewtonMethod} Let $X\subset\mathbb R^n$
be a convex, compact set, $f\colon N\to\mathbb R^n$ be smooth and
fix a point $x\in N$. Let us denote by
\begin{equation}\label{eq:defNewtonOperator}
N(f,X,x) = x - \left[Df(X)\right]_I^{-1}f(x)
\end{equation}
the Interval Newton Operator for a map $f$ on set $X$ with fixed
$x\in X$. Then
\begin{itemize}
\item if $N(f,X,x)\subset\inte X$ then the map $f$ has unique zero
in $X$. Moreover, if $x_*$ is such unique zero of $f$ in $X$ then
$x_*\in N(f,X,x)$.
\item if $N(f,X,x)\cap X=\emptyset$ then the map $f$ has no zeros in
$X$.
\end{itemize}
\end{theorem}

\begin{theorem}(\textbf{Interval Krawczyk Method} \cite{A,Kr,N})
\label{thm:KrawczykMethod}
Let $X\subset\mathbb R^n$ be a convex, compact set, $f\colon
N\to\mathbb R^n$ be smooth and fix a point $x\in N$. Let
$C\in\mathbb R^{n\times n}$ be an isomorphism. Let us denote by
\begin{equation}\label{eq:defKrawczykOperator}
K(f,C,X,x) = x - Cf(x)+(\Id- C\cdot [Df(X)]_I)(X-x)
\end{equation}
the Interval Krawczyk Operator for a map $f$ on set $X$ with fixed
$x\in X$ and matrix $C$. Then
\begin{itemize}
\item if $K(f,C,X,x)\subset\inte X$ then the map $f$ has unique zero
in $X$. Moreover, if $x_*$ is such unique zero of $f$ in $X$ then
$x_*\in K(f,C,X,x)$.
\item if $K(f,C,X,x)\cap X=\emptyset$ then the map $f$ has no zeros in
$X$.
\end{itemize}
\end{theorem}

\subsection{Continuation close to the bifurcation point}

\begin{lemma}\label{lem:c3continuation}
Assume $f_\nu\colon \mathbb R\times\mathbb R^{n-1} \supset X\times Y
\to\mathbb R^n$, $\nu\in[\nu_1,\nu_2]$ is $C^k$ function
both with respect to argument and parameter, with $k \geq
3$, such that
\begin{enumerate}
\item for $\nu\in[\nu_1,\nu_2]$ there exists unique fixed point
$(x_{fp}(\nu),y_{fp}(\nu))$ for $f_\nu$ in $X\times Y$

\item for all $(\nu,x)\in[\nu_1,\nu_2]\times X$ there exists unique
$y(\nu,x)\in \inte Y$ solving equation $y-\pi_y(f^2_\nu(x,y))=0$
and the map $y\colon[\nu_1,\nu_2]\times X\to Y$ is of class
$C^{k}$.

\item the map $G(\nu,x)=x-\pi_x(f^2_\nu(x,y(\nu,x)))$ satisfies
\begin{eqnarray}
    \frac{\partial^3 G}{\partial x^3}(\nu,x)>0,&
    \text{for }\nu\in[\nu_1,\nu_2], x\in X\label{eq:d3g}\\
    G(\nu,x_1)<0,\ G(\nu,x_2)>0,&
    \text{for }\nu\in[\nu_1,\nu_2]\label{eq:opositeSigns}\\
    \forall \nu\in[\nu_1,\nu_2]\
        \exists x_-\in(x_{fp}(\nu),x_2) \label{eq:xMinus}&
        G(\nu,x_-)<0\\
    \forall \nu\in[\nu_1,\nu_2]\
        \exists x_+\in(x_1,x_{fp}(\nu)) &
        G(\nu,x_+)>0\label{eq:xPlus}
\end{eqnarray}
\end{enumerate}
Then there exist two $C^k$ curves $c_1,c_2\colon
[\nu_1,\nu_2]\to\mathbb R^n$ such that for $\nu\in[\nu_1,\nu_2]$
holds $\pi_x(c_1(\nu))<x_{fp}(\nu)<\pi_x(c_2(\nu))$ and $c_i(\nu)$
is a period two point for $f_\nu$, $i=1,2$.

Moreover, if for some $\nu_0\in \{ \nu_1,\nu_2 \}$ holds
$f_{\nu_0}(c_1(\nu_0))=c_2(\nu_0)$ or
$f_{\nu_0}(c_2(\nu_0))=c_1(\nu_0)$ then for all
$\nu\in[\nu_1,\nu_2]$
\begin{equation}\label{eq:c1c2sameorbit-p}
f_{\nu}(c_1(\nu))=c_2(\nu),\quad f_{\nu}(c_2(\nu))=c_1(\nu)
\end{equation}

\end{lemma}

\noindent
\textbf{Proof:} The second assumption and (\ref{eq:d3g})
imply that for a fixed $\nu$ the map $f^2_\nu$ has at most three
fixed points in $X\times Y$. From the first assumption we know
that $f_\nu$ has unique fixed point $(x_{fp}(\nu),y_{fp}(\nu))$ in
$X\times Y$. Therefore any zero of $G(\nu,\cdot)$, which is
different from $(\nu,x_{fp}(\nu))$ corresponds to a period two
point of $f_\nu$. From the continuity of $G$ and from
(\ref{eq:opositeSigns}--\ref{eq:xPlus}) it follows that $G$ has
one zero in each of  the intervals
$x^1_{per}(\nu)\in(x_1,x_{fp}(\nu))$ and
$x^2_{per}(\nu)\in(x_{fp}(\nu),x_2)$. It is easy to see that
functions $x^i_{per}$ are continuous for $i=1,2$. We set
$c_i(\nu)=(x^i_{per}(\nu),y(\nu,x^i_{per}(\nu))$.

We will show the smoothness of $x^i_{per}$, which together with
assumption that $y(x,\nu)$ is $C^k$ implies the smoothness of
$c_i$. It is enough to show that
\begin{displaymath}
\frac{\partial G}{\partial x}(\nu,x^i_{per}(\nu)) \neq 0,
\end{displaymath}
 because
then we can apply the implicit function theorem to obtain the
required differentiability. Let us fix $\nu\in [\nu_1,\nu_2]$.
Observe that condition (\ref{eq:d3g}) implies that for any fixed
$\nu$ the function $x \mapsto \frac{\partial G }{\partial
x}(\nu,x)$ has at most two zeros in $[x_1,x_2]$.  From remaining
assumptions it is clear that on interval $[x_1,x_{fp}(\nu)]$ and
$[x_{fp}(\nu),x_2]$ the function $x \mapsto G(\nu,x)$ has strictly
positive maximum and strictly negative minimum, respectively.
Therefore these extremal points are zeros of $\frac{\partial
G}{\partial x}(\nu,x)$ and obviously they are different from
points $x^i_{per}(\nu)$, which are zeros of $G(\nu,\cdot)$. Hence
we have shown that $\frac{\partial G}{\partial
x}(\nu,x^i_{per}(\nu))\neq 0$.

Assume that $\nu_0=\nu_1$ (the other case is analogous). From the
implicit function theorem it follows that for some $\nu' > \nu_1$
condition (\ref{eq:c1c2sameorbit-p}) is satisfied for $\nu_1 \leq
\nu < \nu'$. Let $\nu_m$ be supremum of such $\nu' \leq \nu_2$. It
is easy to see that $\nu_m=\nu_2$, because at $\nu_m$ is also
satisfied by the continuity and implicit function theorem allows us
to extend the range of $\nu$ satisfying (\ref{eq:c1c2sameorbit-p})
to the right if $\nu_m < \nu_2$.

 \qed

\section{Extracting the dynamical information from Lia\-pu\-nov-Schmidt reduction}
\label{sec:dyninfo}

As was mentioned already in  Section~\ref{subsec:hypbif} the
Liapunov-Schmidt projection does not give us any direct
information about the dynamics of bifurcating solutions regarding
the invariant manifolds of the bifurcating objects as required by
Def.~\ref{def:perdoubl}.
 In this section
following the ideas of de Oliveira and  Hale \cite{H,OH}, we show that the information
obtained from the Liapunov-Schmidt reduction and  the
spectrum of the bifurcating fixed point curve is enough  to say
precisely, what is the dynamics in the neighbourhood of the
bifurcation point.

Our argument  follow the ideas from \cite[Chapter 9, Thm. 3.1 and
4.2]{CH}, where an analogous problem was considered for fixed points
for ODEs and periodic orbits for periodically forced ODEs. The
notion of the central manifold \cite{K} plays crucial role in this
proof.

\begin{figure}
\centerline{\includegraphics{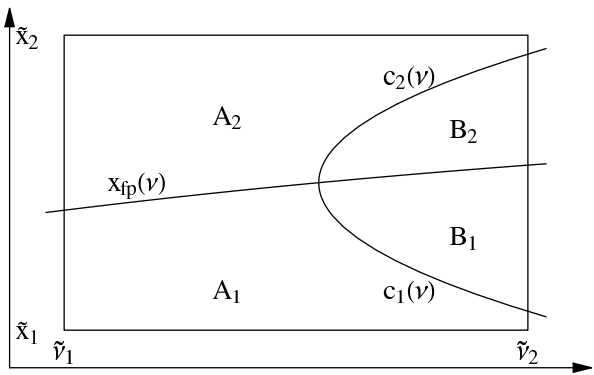}
    }
    \caption{Location of sets $A_i$, $B_i$ with respect to zeros of $G$.\label{fig:bifsets}}
\end{figure}

\begin{theorem}\label{thm:dynInfo}
Let  $P_\nu:\mathbb{R}^n \to \mathbb{R}^n$ for $\nu \in
[\nu_1,\nu_2]$ be a $C^k$-map ($k \geq 3$) both with respect to
$\nu$ and its arguments. Assume that on the set $V=[\nu_1,\nu_2]
\times [x_1,x_2] \times Y$, where $Y \subset \mathbb{R}^{n-1}$ is
a closure of an open set, we were able to perform the
Liapunov-Schmidt reduction and verify assumptions of
  Theorem~\ref{thm:per-doubling}. Let $(\nu_b,z_{b})$ be the bifurcation point
  and $z_{fp}(\nu)=(x_{fp}(\nu),y(\nu,x_{fp}(\nu))$ be the fixed
  point curve for $P_\nu$ in $V$.

Let $v$ be the eigenvector of $\frac{\partial P_{\nu_b}}{\partial
z}(z_b)$ corresponding to the eigenvalue $-1$. We assume that
$\pi_x v \neq 0$.

 Assume that there exists $\epsilon >0$, such that
for all $\nu \in [\nu_1,\nu_2]$ holds
\begin{equation*}
 \mbox{Sp}\left(\frac{\partial P_{\nu}}{\partial
 z}(z_{fp}(\nu))\right) =   A_\nu \cup B_\nu \cup \{ \lambda(\nu) \},
\end{equation*}
 where $\lambda(\nu) \in \mathbb{R}$   has the multiplicity one,
 $A_\nu \subset \{\alpha \in \mathbb{C},\quad |\alpha| < 1 - \epsilon
\}$
  and $B_\nu \subset  \{\beta \in \mathbb{C},\quad |\beta| > 1 + \epsilon \}$.
  Moreover, we assume
 that
 \begin{eqnarray*}
   \lambda(\nu_1) &\subset& (-1, 1) \\
   \lambda(\nu_2) &<& -1.
 \end{eqnarray*}

Then the map $P$ has a period doubling bifurcation at
$(\nu_b,x_b,y(\nu_b,x_b))$.
\end{theorem}
\textbf{Proof:} Let $\nu:[\bar{x}_1,\bar{x}_2] \to [\nu_1,\nu_2]$
be the function from assumption A3 of
Theorem~\ref{thm:per-doubling} (in fact of Lemma~\ref{lem:Gsol})
is satisfied. In the notation used in
Theorem~\ref{thm:per-doubling} we have
$(\nu_b,z_{b})=(\nu(\bar{x}_0),(\bar{x}_0,y(\nu(\bar{x}_0),\bar{x}_0)))$.
 Let $c_1(\nu)$ and $c_2(\nu)$ be respectively lower and upper
branch of the graph of the function $x \to \nu(x)$ giving period-2
points -- see Fig.~\ref{fig:bifsets}.

Let us define a map $H:V \to \mathbb{R} \times \mathbb{R}^{n}$ by
\begin{equation*}
  H(\nu,z)=(\nu,P(\nu,z)).
\end{equation*}
Consider the spectrum of $DH(\nu_b,z_{b})$. It is easy to see that
$+1$ is an eigenvalue of  $DH(\nu_b,z_{b})$ of multiplicity one,
$\lambda=-1$ has also multiplicity one and all other eigenvalues are
off the unit circle.

We apply the center manifold theorem \cite{K,G,HPS} to $H$ in the
neighbourhood of $(\nu_b,z_b)$. Therefore, there exists a
neighbourhood $M$ of $(\nu_b,z_b)$ and  two-dimensional center
manifold $W^c \subset M$ such that
\begin{eqnarray*}
  \forall (\nu,z) \in W^c \quad \mbox{if $H^i(\nu,z) \in M$ then $H(\nu,z) \in W_c$, for
  $i=-1,1$}\\
  \Inv(M,H) \subset W^c
\end{eqnarray*}
$W^c$ is tangent at $(\nu_b,z_b)$ to the subspace spanned by
vectors $\{(1,0),(0,v)\} \subset \mathbb{R} \times \mathbb{R}^n$.
Observe that from our assumption about $v$, i.e. $\pi_x(v)\neq0$,
it follows that we can use on $W^c$ in the neighbourhood of
$(\nu_b,z_{b})$ the same coordinates $(\nu,x)$ as in the
Liapunov-Schmidt reduction. There exists a neighbourhood of
$(\nu_b,z_{b})$ denoted by $U=[\tilde{\nu}_1,\tilde{\nu_2}] \times
[\tilde{x}_1,\tilde{x}_2] \times \tilde{Y} \subset M \cap V$ and
$C^k$-functions $h:[\tilde{\nu}_1,\tilde{\nu_2}] \times
[\tilde{x}_1,\tilde{x}_2] \to  \tilde{Y}$ and
$f:[\tilde{\nu}_1,\tilde{\nu_2}] \times [\tilde{x}_1,\tilde{x}_2]
\to \mathbb{R}$ satisfying
\begin{eqnarray}
  W^c&=&\{ (\nu,x,h(\nu,x)) \} \nonumber\\
  P(\nu,x,h(\nu,x))&=& (f(\nu,x),h(\nu,f(\nu,x))) \nonumber\\
  \Inv(U,H) &\subset& W^c.  \label{eq:InvinWc}
\end{eqnarray}

Let us stress  that the dynamics of $P_\nu$ in $W^c$  is
one-dimensional , namely that of $x \mapsto f(\nu,x)$.

From the Liapunov-Schmidt reduction we know that a point $(\nu,x,y)
\in U$ has period one or two with respect to map $H$ iff
$y=y(\nu,x)$ and
 $G(\nu,x)=0$. Let $N=[\tilde{\nu}_1,\tilde{\nu_2}] \times
[\tilde{x}_1,\tilde{x}_2]$. If $U$ is chosen to be sufficiently
close to the bifurcation point, then the set $N \setminus \{(\nu,x)
\: | \: G(\nu,x)=0 \}$ has four connected components -- see
Fig~\ref{fig:bifsets}. Namely,
\begin{eqnarray*}
  A_1 &=& \left\{ (\nu,x) \in N \, | \, \left( (\nu \leq \nu_b) \mbox{ and } (x<x_{fp}(\nu))  \right)  \mbox{ or } \right. \\
     & & \left. \left(
  (\nu > \nu_b) \mbox{ and } (x < c_1(\nu))  \right)  \right\}, \\
  A_2 &=& \left\{ (\nu,x) \in N \, | \, \left( (\nu \leq \nu_b) \mbox{ and } (x>x_{fp}(\nu))  \right)  \mbox{ or } \right. \\
      & &    \left.  \left(
  (\nu > \nu_b) \mbox{ and } (x > c_2(\nu))  \right)  \right\}, \\
  B_1 &=& \left\{ (\nu,x) \in N \, | \,  (\nu > \nu_b) \mbox{ and } ( x_{fp}(\nu)> x > c_1(\nu))   \right\}, \\
  B_2 &=& \left\{ (\nu,x) \in N \, | \,  (\nu > \nu_b) \mbox{ and } ( x_{fp}(\nu) < x < c_2(\nu))
  \right\}.
\end{eqnarray*}

We also require that
\begin{equation}
  \tilde{x}_1 < c_1(\tilde{\nu}_2) < c_2(\tilde{\nu}_2) <
  \tilde{x}_2. \label{eq:shapeU}
\end{equation}

On each of these components the function
$d(\nu,x)=x-f(\nu,f(\nu,x))$ must have a constant sign. Observe
that on $A_2$  we have
\begin{equation}
  x - f(\nu,f(\nu,x)) > 0, \quad \mbox{for $(\nu,x) \in A_2$}
  \label{eq:sgnA2}
\end{equation}
because $z_{fp}(\tilde{\nu}_1)$ is attracting on $W^c$ and we
consider the second iterate. Analogously we obtain
\begin{equation}
  x - f(\nu,f(\nu,x)) < 0, \quad \mbox{for $(\nu,x) \in A_1$}
  \label{eq:sgnA1}
\end{equation}

For the component $B_2$  we have
\begin{equation*}
  x - f(\nu,f(\nu,x)) < 0, \quad \mbox{for $(\nu,x) \in B_2$}
\end{equation*}
because $x_{fp}(\bar{\nu}_2)$ is repelling on $W^c$ and we consider
the second iterate. Analogously
\begin{equation*}
  x - f(\nu,f(\nu,x)) > 0, \quad \mbox{for $(\nu,x) \in B_1$}.
\end{equation*}

For a subset $Z\subset N$ by $Z_\nu$ we will denote $Z_\nu=\{x:
(\nu,x)\in N\}$. Observe that for each $\nu \in
[\tilde{\nu}_1,\tilde{\nu}_2]$ holds
\begin{equation}
  f_\nu^2((A_i)_\nu) \cap [\tilde{x}_1,\tilde{x}_2] \subset (A_i)_\nu, \quad
           f_\nu^2((B_i)_\nu) \cap  [\tilde{x}_1,\tilde{x}_2]  \subset
           (B_i)_\nu
    \quad i=1,2. \label{eq:invcomponets}
\end{equation}

For the proof of (\ref{eq:invcomponets}) observe that map
$l(\nu,x)=(\nu,f(\nu,x))$ maps connected components of $N \setminus
\{G(\nu,x)=0\}$ into connected components, i.e. for any $S \in
\{A_1,A_2,B_1,B_2\}$ there exists $T=T(S) \in \{A_1,A_2,B_1,B_2\}$,
such that
\begin{equation}
\left(l(S) \cap N \right) \subset T(S), \label{eq:cc-into-cc}
\end{equation}
because
\begin{equation*}
l(G^{-1}(0)\cap N)\cap N = G^{-1}(0)\cap N = l^{-1}(G^{-1}(0)\cap
N)\cap N
\end{equation*}
Observe that the relevant eigenvalue of $DP_\nu$ at $z_{fp}(\nu)$
describing the dynamics on $W^c$ is $\lambda(\nu)$, which is real
and since we consider the second iterate we see that in the
neighbourhood of the fixed point curve we have points mapped into
the same component. This  together with (\ref{eq:cc-into-cc}) proves
(\ref{eq:invcomponets}).

From the above considerations we obtain for $\nu \leq \nu_b$
\begin{eqnarray*}
   x_{fp}(\nu) < f(\nu,f(\nu,x)) < x , \quad \mbox{for $   x_{fp}(\nu) <  x  \leq
   \tilde{x}_2$} \\
   x_{fp}(\nu) > f(\nu,f(\nu,x)) > x , \quad \mbox{for $   x_{fp}\nu) >  x
   \geq
   \tilde{x}_1$}.
\end{eqnarray*}
The above conditions, (\ref{eq:InvinWc}) and nonexistence of other
period two points in $U$ imply that
\begin{equation*}
  \Inv(\pi_{x,y}U,P_\nu)=\{z_{fp}(\nu)\}, \quad \mbox{for $\nu \in
  [\tilde{\nu},\nu_b]$}.
\end{equation*}

For $\nu \in (\nu_b,\tilde{\nu}_2]$ we have
\begin{eqnarray*}
   x_{fp}(\nu) < x < f(\nu,f(\nu,x)) < c_2(\nu)   , \quad \mbox{for $   x_{fp}(\nu) <  x
   < c_2(\nu) $ } \\
  x_{fp}(\nu) > x > f(\nu,f(\nu,x)) > c_1(\nu)   , \quad \mbox{for $   x_{fp}(\nu) >  x
   > c_1(\nu) $ } \\
   c_{2}(\nu) < f(\nu,f(\nu,x)) < x   , \quad \mbox{for $   x   > c_2(\nu) $ } \\
   c_{1}(\nu) > f(\nu,f(\nu,x)) > x   , \quad \mbox{for $   x   < c_1(\nu) $
   }.
\end{eqnarray*}
The above conditions, conditions
(\ref{eq:shapeU},\ref{eq:InvinWc}) and nonexistence of other
period two points in $U$ imply that for $\nu \in
(\nu_b,\tilde{\nu}_2]$
\begin{eqnarray*}
 \Inv(\pi_{x,y}U,P_\nu)=\{ (x,h(\nu,x)) \, | \, x \in [c_1(\nu),c_2(\nu)] \}
\end{eqnarray*}

 \qed

 We would like to stress here that, contrary to all
previous theorems and lemmas, in the proof of the above theorem we
prove the statements about the invariant manifold of bifurcating
orbits from Definition~\ref{def:perdoubl} on some set $U$, whose
size we do not  control, whereas it is given by the range of the
existence of the central manifold. In principle, this range can be
inferred from the proof of the center manifold theorem, but it  will
be an interesting task to develop a computable approach, which will
allow to rigorously prove these facts  on the whole set $V$. Such
task will require explicite estimates about the central manifold in
the region very close to the bifurcation and some other tools, may
be of Conley index type \cite{MM}, further away from the
bifurcation.

\section{Application to the R\"ossler system.}\label{sec:firstBif}
Consider an autonomous ODE in $\mathbb R^3$ called the R\"ossler
system \cite{R}
\begin{equation}\label{eq:rosslerSystem}
\left\{
    \begin{array}{rcl}
        x' &=& -y-z\\
        y' &=& x+by\\
        z' &=& b+z(x-a)
    \end{array}
\right.
\end{equation}
The classical parameter values (considered by R\"ossler) are $a=5.7$
and $b=0.2$. For the remainder of this paper we fix $b=0.2$.

\begin{figure}[htbp]
    \centerline{\includegraphics{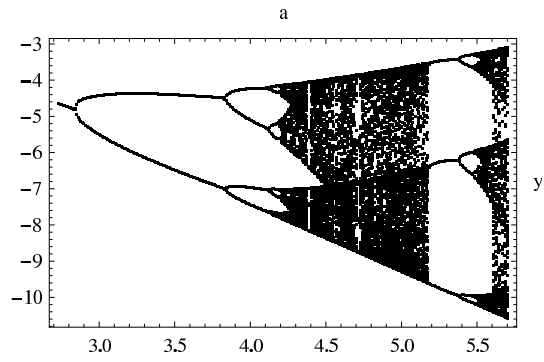}}
    \caption{Bifurcation diagram for the R\"ossler system\label{fig:bifDiagram}}
\end{figure}

\begin{figure}[htbp]
    \centerline{
        \begin{tabular}{cc}
        \includegraphics[width=.5\textwidth]{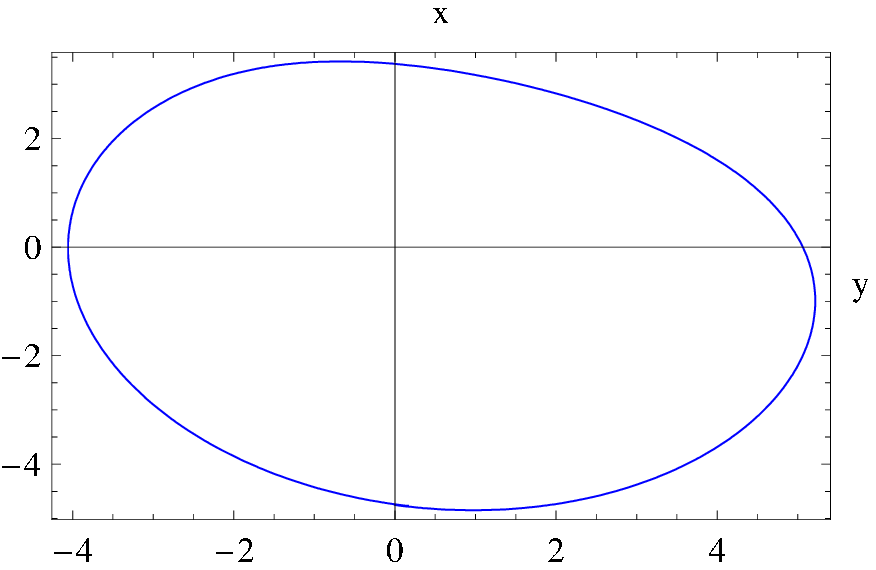}&
        \includegraphics[width=.5\textwidth]{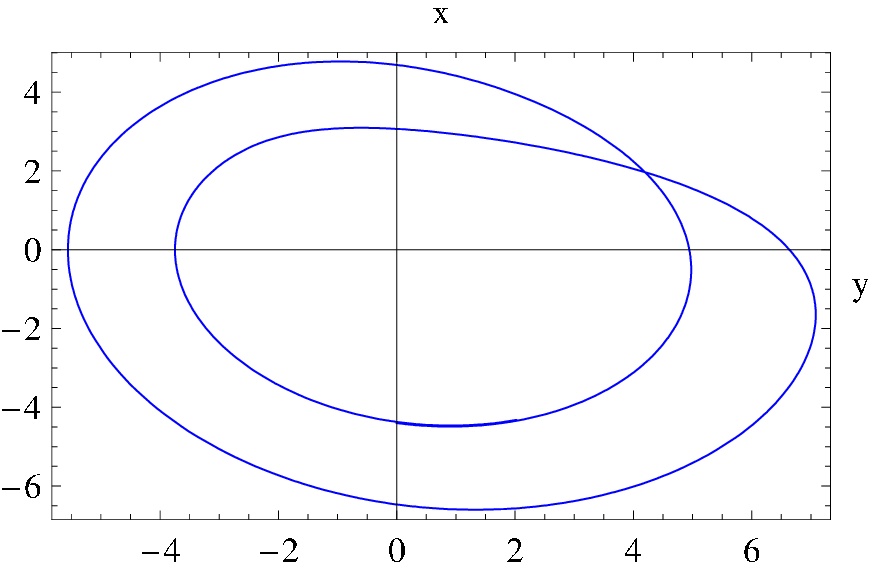}\\
        \includegraphics[width=.5\textwidth]{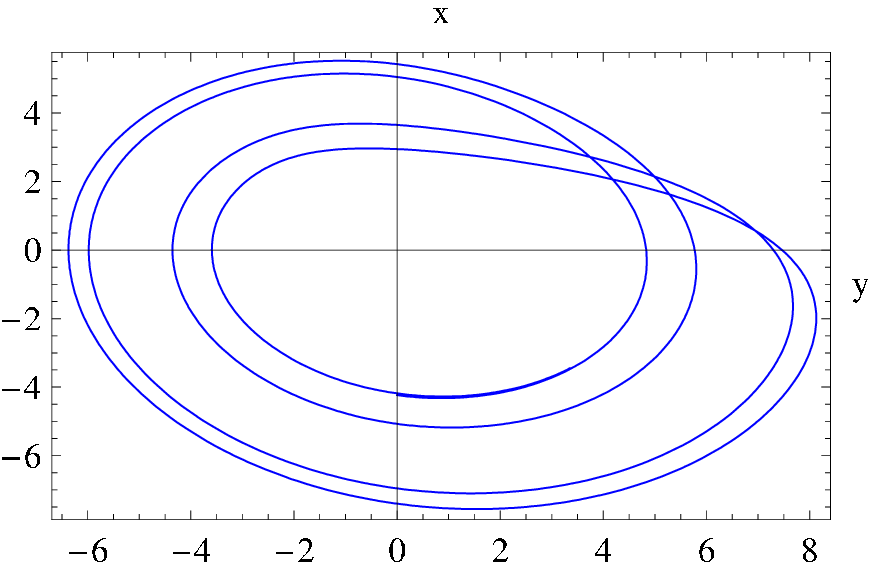}&
        \includegraphics[width=.5\textwidth]{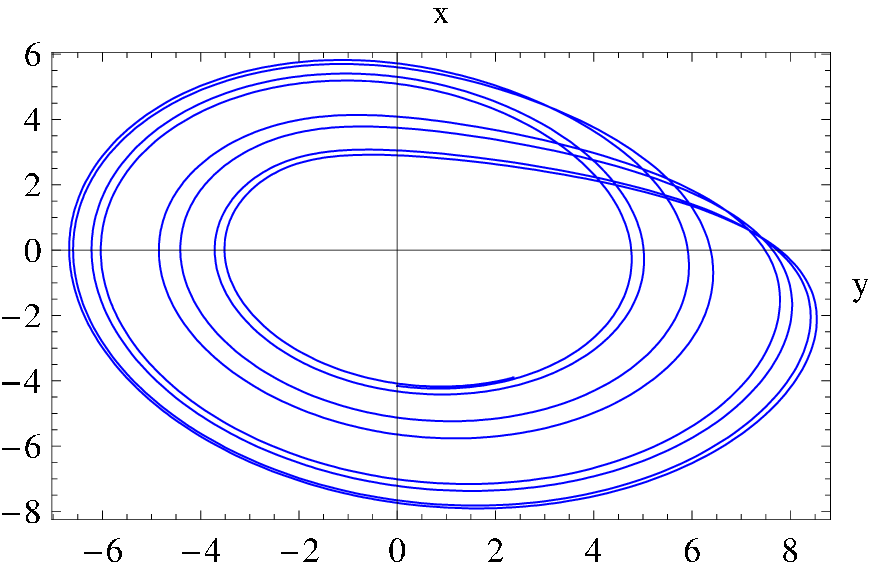}
        \end{tabular}
        }
    \caption{Periodic orbits corresponding to fixed point, period two point, period four point
    and period eight point for the Poincar\'e map.
    Parameter values are $a=2.8$, $a=3.5$, $a=4$ and $a=4.2$.\label{fig:perOrbits}}
\end{figure}

The system (\ref{eq:rosslerSystem}) has been extensively studied
in the literature numerically and is treated in the literature as
one of classical examples of systems generating chaotic attractor.
Yet, the number of rigorous results concerning it is very small.
In Fig.~\ref{fig:bifDiagram} we show a numerically obtained
bifurcation diagram for periodic orbits on section $x=0$ with
$b=0.2$ and $a$ as parameter. We see that when the parameter $a$
increases from $2$ to $5.7$ one observes a cascade of period
doubling bifurcations. In Fig.~\ref{fig:perOrbits} we show some
periodic orbits for different values of $a$. Our goal in this
section is to validate the part of the bifurcation diagram in
Fig.\ref{fig:bifDiagram} containing two first period doublings
using the approach introduced in the previous sections.

Let us list the few known rigorous results about (\ref{eq:rosslerSystem}).
 Pilarczyk (see \cite{P}
and references given there) gave a computer assisted proof of the
following facts: for $a=2.2$ there exists periodic orbit, for
$a=3.1$ there exists two periodic orbits. However from his proof
one cannot infer any information about the dynamical character of
these orbits. He constructs suitable isolating neighborhoods,
which have an index of an attracting or a hyperbolic orbit with
one unstable direction, but no such claim can be made about the
periodic orbit proved to exists. In fact we do not even known,
wether this orbit is unique.

Finally, for the classical parameter values $b=0.2$ and $a=5.7$ the
system is chaotic \cite{Z3}  in the following sense: a suitable
Poincar\'e map has an invariant set $S$ and the dynamics on $S$  contains
the shift map  dynamics on three symbols.

Before proceeding any further  we need to introduce some notation.
Let $\Pi=\{(x,y,z)\in\mathbb R^3\ |\ x=0,x'>0\}$ be a Poincar\'e
section. Since for $u\in\Pi$ the first coordinate is equal to zero
we will use the remaining two coordinates $u=(y,z)$ to represent a
point on $\Pi$. For a fixed parameter value $a>0$ by
$P_a=(P_{a,y},P_{a,z})\colon \Pi\to\Pi$ we will denote the
corresponding Poincar\'e return map. By $P$ we will denote the map
defined by $P(a,y,z)=(a,P_{a,y}(y,z),P_{a,z}(y,z))$.

Apparently the first period doubling bifurcation is observed for
$a\approx2.832445$ and the second one for $a\approx 3.837358$. In
the remainder of this section we discuss the computer assisted
proof of the existence of both these bifurcations. In our
presentation we will discuss the first one more in details, while
for the second one we will just state relevant lemmas and
estimates.

Let $u_0=(y_0,z_0)$ be an approximate fixed point for $P_{a_0}$,
i.e. we set
\begin{equation}\label{eq:firstBifPoint}
u_0=(y_0,z_0) = (-4.7946653021070986256,0.052488098609082899093)
\end{equation}
and put
\begin{equation}\label{eq:firsBifCoordSystem}
M = \begin{bmatrix}
    0.99999765967819775891 &-0.9582095926217468751\\
    0.0021634782474835700244 &-0.28606708410382636343
\end{bmatrix}
\end{equation}
The columns of $M$ are normalized approximate eigenvectors of
$DP_{a_0}(u_0)$, where first column corresponds to the eigenvalue
close to $-1$ and the second one to the eigenvalue close to zero.
On section $\Pi$ we choose new coordinates $(\tilde y,\tilde z) =
M^{-1}((y,z)-u_0)$ and since, in the sequel, we will use only the
new coordinates we will drop the tilde.

Define
\begin{equation*}
\begin{array}{rcl}
     A &=& [a_1,a_2]=[2.83244,2.832446]\\
     Y &=& [y_1,y_2]=1.3107\cdot[-1,1]\cdot 10^{-3}\\
     Z &=& [z_1,z_2]=1.3107\cdot[-1,1]\cdot 10^{-4}
\end{array}
\end{equation*}

Now our goal is present the proof of  the following theorem
\begin{theorem}\label{thm:rosslerFirstPD}
The map $P_a$ has a period doubling bifurcation at some point
$(a,y,z) \in \inte ( A\times Y\times Z)$.
\end{theorem}

\begin{rem}
The existence of period doubling bifurcation is a local phenomenon.
In fact the sets $A$, $Y$, $Z$ can be chosen to be smaller which
speed up the proof ($13$ minutes versus $87$ minutes), namely we
were able to prove the existence of period doubling bifurcation in
the set
\begin{equation*}
\begin{array}{rcl}
     A &=& [a_1,a_2]=[2.83244,2.832445028]\\
     Y &=& [y_1,y_2]=[-1,1]\cdot 10^{-4}\\
     Z &=& [z_1,z_2]=[-1,1]\cdot 10^{-5}
\end{array}
\end{equation*}
However, the choice of larger set facilitates the  proof of the
existence of connecting branch of period two points between first
and second period doubling bifurcation, because decreasing $a_2$
results in the eigenvalue of period-two points to be very close to
$1$, which makes it very difficult to rigorously continue it.
\end{rem}

\subsection{The existence of fixed point curve.}
\begin{lemma}\label{lem:fixedPointCurve}
There exists function $(y_{fp},z_{fp})\colon A\to Y\times Z$
of class $C^\infty$ such that for $(a,y,z)\in A\times Y\times Z$
holds
\begin{equation*}
P_a(y,z)=(y,z) \qquad \text{iff}\qquad (y,z)=(y_{fp}(a),z_{fp}(a))
\end{equation*}
and
\begin{equation}\label{eq:yfpPrimeBound}
y_{fp}'(A)\subset[-1.3336825610133946629,-1.3275439332565022177]
\end{equation}
\end{lemma}
\textbf{Proof:} The proof, which is computer assisted, consists from
two parts, in the first one we prove the existence of the fixed
point curve and in the second part we establish estimate
(\ref{eq:yfpPrimeBound}).

For the first part, we use the Interval Newton Method
(Theorem~\ref{thm:NewtonMethod}) and $C^1$-Lohner algorithm to
prove  that for $a\in A$ there exists a unique fixed point
$(y_{fp}(a),z_{fp}(a))$ for $P_a$ in $Y\times Z$. In computations
we insert the whole set $A\times Y\times Z$ as an initial
condition in our routine, which computes the Interval Newton
Operator and obtain that the for all $a\in A$ the fixed point
$(y_{fp}(a),z_{fp}(a))$ belongs to the set
\begin{multline}\label{eq:newtonFixedPointCurve}
N:=N(\Id-P_a,Y\times Z)=\\\begin{bmatrix}
[-2.838378938597049559,3.2727784971172813446]\cdot 10^{-5}\\
[-4.8121450471307824034,4.2979575521536656939]\cdot
10^{-6}\end{bmatrix}^T
\end{multline}

To obtain (\ref{eq:yfpPrimeBound})  we apply $C^1$-Lohner
algorithm \cite{Z1} to the system
\begin{equation}\label{eq:extendedRosslerSystem}
\left\{
\begin{array}{rcl}
        x' &=& -y-z\\
        y' &=& x+by\\
        z' &=& b+z(x-a)\\
        a' &=& 0
\end{array}
\right.
\end{equation}
with $b=0.2$ in order to compute a bound for $y_{fp}'(A)$.
Differentiating
\begin{equation*}
P(a,y_{fp}(a),z_{fp}(a))=(a,y_{fp}(a),z_{fp}(a))
\end{equation*}
with respect to $a$ we obtain
\begin{equation}\label{eq:yfpPrimeFormula}
    y_{fp}' = \frac{\frac{\partial P_{a,y}}{\partial a}\left(
        1-
        \frac{\partial P_{a,z}}{\partial z}\right)
        + \frac{\partial P_{a,z}}{\partial a} \cdot  \frac{\partial P_{a,y}}{\partial z}
    }
    {
        \left(\frac{\partial P_{a,y}}{\partial y}-1\right)
        \left(\frac{\partial P_{a,z}}{\partial z}-1\right)
        -\frac{\partial P_{a,y}}{\partial z}
        \frac{\partial P_{a,z}}{\partial y}
    }
\end{equation}
where the partial derivatives of $P$ are evaluated at
$(y_{fp}(a),z_{fp}(a))$.

We use the set $A\times N$, where $N$ is defined in
(\ref{eq:newtonFixedPointCurve}) as initial condition in our routine
which computes partial derivatives of $P$ and after substituting
them to (\ref{eq:yfpPrimeFormula}) we obtain a bound for $y_{fp}'$
as in (\ref{eq:yfpPrimeBound}).

We used the Taylor method of order $14$ and the time step equal to
$0.02$ to integrate the system (\ref{eq:rosslerSystem}) in $\mathbb
R^3$ for the first part of the proof and the order $10$ and the time
step $0.01$ when we integrate the extended system
(\ref{eq:extendedRosslerSystem}) in the second part. \qed

\begin{lemma}\label{lem:derEigenvalue}
The eigenvalues  $\lambda_1,\lambda_2:A\to\mathbb R$   of
$DP_a(y_{fp}(a),z_{fp}(a))$ are given by
\begin{eqnarray*}
    \lambda_1(a) &=& \frac{1}{2}\left(
        \frac{\partial P_{a,y}}{\partial y} + \frac{\partial P_{a,z}}{\partial
        z}-s(a)
    \right),\\
    \lambda_2(a) &=& \frac{1}{2}\left(
        \frac{\partial P_{a,y}}{\partial y} + \frac{\partial P_{a,z}}{\partial
        z}+s(a)
    \right),\\
    s(a) &=& \sqrt{\left(\frac{\partial P_{a,y}}{\partial y}-\frac{\partial P_{a,z}}{\partial z}\right)^2+4\frac{\partial P_{a,y}}{\partial z}\frac{\partial P_{a,z}}{\partial y}}
\end{eqnarray*}
where partial derivatives of $P$ are evaluated at
$(y_{fp}(a),z_{fp}(a))$. Let $v(a)$ be the normalized eigenvector
corresponding to eigenvalue $\lambda_1(a)$. Then
\begin{eqnarray*}
    \lambda_1(a_1)&\in&   [-0.99999781944914578613,-0.99999548919217751131]\\
    \lambda_1(a_2)&\in&   [-1.00000064581599335217,-1.00000064581598072628] \\
    \lambda_2(A)&\subset& [-0.0013533261367103342071,0.0013530378340487671934]\\
    \lambda_1'(A)&\subset&[-0.70107900728585614836,-0.62770519734197127715]\\
    v_y(A)&\subset& \pm[0.99728887963031764841,1.0027184248801992439]
\end{eqnarray*}
where $v_y$ denotes the $y$ coordinate of $v$.
\end{lemma}
\textbf{Proof:} We leave the derivation of formulas for
$\lambda_1,\lambda_2$  to the reader.
 We used the $C^1$-Lohner algorithm applied to the
system (\ref{eq:rosslerSystem}) in order to compute bounds for
$\lambda_1(a_1)$ and $\lambda_1(a_2)$. Since the parameter $a_2$ has
been chosen to be very close to the bifurcation parameter we find
difficulties with the verification of condition $\lambda_1(a_2)<-1$ in
computations performed in interval arithmetics based on double
precision (52-bit mantissa) boundary value type. In our computations
we used interval arithmetics based on float numbers with $150$-bit
mantissa (MPFR \cite{MPFR} and GMP \cite{GMP} packages).

Since the eigenvalue $\lambda_1(a)$ of $DP_a(y_{fp}(a),z_{fp}(a))$
is given by an explicit formula one can  express $\lambda_1'(a)$ in
terms of first and second order partial derivatives of $P$. We
obtain
\begin{eqnarray*}
\lambda_1'(a) &=& \frac{1}{2}\left(\frac{\partial}{\partial
a}\frac{\partial P_{a,y}}{\partial y}+\frac{\partial}{\partial
a}\frac{\partial P_{a,z}}{\partial z}-s'(a)\right)\\
s'(a) &=& \frac{1}{s(a)}\left(\frac{\partial P_{a,y}}{\partial
y}-\frac{\partial P_{a,z}}{\partial
z}\right)\left(\frac{\partial}{\partial a}\frac{\partial
P_{a,y}}{\partial y}-\frac{\partial}{\partial a}\frac{\partial
P_{a,z}}{\partial z}\right)\\
&+&\frac{2}{s(a)}\left(
    \left(\frac{\partial}{\partial a}\frac{\partial P_{a,y}}{\partial
    z}\right)\frac{\partial P_{a,z}}{\partial y}
    +
    \left(\frac{\partial}{\partial a}\frac{\partial P_{a,z}}{\partial
    y}\right)\frac{\partial P_{a,y}}{\partial z}
\right)
\end{eqnarray*}
where the symbols $\frac{\partial}{\partial a}\frac{\partial
P_{a,y}}{\partial z}$ and $\frac{\partial}{\partial a}\frac{\partial
P_{a,z}}{\partial y}$ should be understood as
\begin{eqnarray*}
\frac{\partial}{\partial a}\frac{\partial P_{a,z}}{\partial
y}(y_{fp}(a),z_{fp}(a))
    &=& \frac{\partial^2 P_{a,z}}{\partial a\partial
    y}(y_{fp}(a),z_{fp}(a))\\
    &+&
    \frac{\partial^2 P_{a,z}}{\partial
    y^2}(y_{fp}(a),z_{fp}(a))y_{fp}'(a)\\
    &+&
    \frac{\partial^2 P_{a,z}}{\partial y\partial
    z}(y_{fp}(a),z_{fp}(a))z_{fp}'(a)
\end{eqnarray*}
and $y_{fp}'(a)$ and $z_{fp}'(a)$ can be computed as in
(\ref{eq:yfpPrimeFormula}). Next, we applied the $C^2$-Lohner
algorithm \cite{WZn} to the extended system
(\ref{eq:extendedRosslerSystem}) in order to compute a bound for
the first and the second order partial derivatives of $P$ and in
consequence a bound for $\lambda_1'(A)$.

We inserted $A\times N$, where $N$ is defined in
(\ref{eq:newtonFixedPointCurve}), as the initial condition in our
routine, which computes the  partial derivatives of Poincar\'e map
up to second order. In these computations we simultaneously
computed bounds for $\lambda_2(A)$ and $\lambda_1'(A)$. The
parameter settings of the Taylor method used in the computations
are listed in Table~\ref{tab:C1C2settings}.

\begin{table}
\begin{center}
\begin{tabular}{|c|c|c|}
    \hline
        \mbox{ } & order & step\\
    \hline
        $\lambda_2(A)$, $\lambda_1'(A)$ & $10$ & $0.03$\\
    \hline
        $\lambda_1(a_1)$ & $10$ & $0.1$\\
    \hline
        $\lambda_1(a_2)$ - 150-bit precision & $14$ & $0.05$\\
    \hline
\end{tabular}
\end{center}
\caption{Parameters of the $C^1-C^2$-Lohner
algorithms.\label{tab:C1C2settings}}
\end{table}
\qed

\subsection{The existence of Liapunov-Schmidt reduction.}
\begin{lemma}\label{lem:LiapSchmidtExistence}
For all $(a,y)\in A\times Y$ there exists unique $z=z(a,y)\in Z$
such that
\begin{equation}
P_{a,z}^2(y,z)=z\qquad \text{iff}\qquad
z=z(a,y)\label{eq:LiapSchmidtExistence}
\end{equation}
and the map $z\colon A\times Y\to Z$ is smooth of class $C^\infty$.
Moreover, the map $G\colon A\times Y\to \mathbb R$  defined by
\begin{equation*}
G(a,y) = y - P^2_{a,y}(y,z(a,y))
\end{equation*}
satisfies
\begin{eqnarray}
\frac{\partial^3 G}{\partial y^3}(A\times Y) \subset  [1.8296823158090675943,7.2204769494502958338]\label{eq:d3gEstimation}\\
\frac{\partial^2 G}{\partial y^2}(A\times Y) \subset
[-0.2084557586786322414,0.2080871792788867581]
\label{eq:d2gEstimation}
\end{eqnarray}
\end{lemma}
\textbf{Proof:}
Let us fix $(a,y)\in A\times Y$ and define a function $V_{a,y}\colon
Z\to \mathbb R$ by $V_{a,y}(z)=z-P_{a,z}^2(y,z)$. The computer
assisted proof of this Lemma consists of the following steps

\begin{itemize}
\item We divide interval $A$ onto $30$ parts. For each
subinterval $\bar A$ in this covering we proceed as follows

\item Using Interval Newton Method (Theorem~\ref{thm:NewtonMethod}) we
verified that for all $(a,y)\in \bar A\times Y$ the function
$V_{a,y}$ has exactly one zero in $Z$. Denote this zero by $z(a,y)$.
This defines the unique map $z\colon \bar A\times Y\to Z$ which is
smooth by implicit function theorem and which satisfies
(\ref{eq:LiapSchmidtExistence}).

\item
Let $\bar Z$ denote a bound for $z(\bar A,Y)$ resulting from the
previous step. Differentiating $z(a,y)-P_{a,z}^2(y,z(a,y))=0$ with
respect to $y$ we obtain
\begin{eqnarray*}
\left(1-\frac{\partial P_{a,z}}{\partial z}\right)\frac{\partial
z}{\partial y} &=&
    \frac{\partial P_{a,z}}{\partial y}\\
\left(1-\frac{\partial P_{a,z}}{\partial z}\right)\frac{\partial^2
z}{\partial y^2} &=&
    \frac{\partial^2 P_{a,z}}{\partial z^2} (\frac{\partial z}{\partial y})^2
    +2 \frac{\partial^2 P_{a,z}}{\partial y\partial z}\frac{\partial z}{\partial y}
    +\frac{\partial^2 P_{a,z}}{\partial y^2}\\
\left(1-\frac{\partial P_{a,z}}{\partial z}\right)\frac{\partial^3
z}{\partial y^3} &=&
    \frac{\partial^3 P_{a,z}}{\partial z^3}\left(\frac{\partial z}{\partial y}\right)^3
    +3 \frac{\partial^3 P_{a,z}}{\partial y\partial z^2} \left(\frac{\partial z}{\partial y}\right)^2
    \\&+&3 \left(
        \frac{\partial^2 z}{\partial y^2}
        \frac{\partial^2 P_{a,z}}{\partial z^2}+\frac{\partial^3 P_{a,z}}{\partial y^2\partial z}
        \right) \frac{\partial z}{\partial y}\\
    &+&3 \frac{\partial^2 z}{\partial y^2}
    \frac{\partial^2 P_{a,z}}{\partial y\partial z}+\frac{\partial^3 P_{a,z}}{\partial y^3}
\end{eqnarray*}
We see that we can compute all the partial derivatives of $z(a,y)$
as a functions of partial derivatives of $P$. Hence, partial
derivatives of $G(a,y) = y-P_{a,y}^2(y,z(a,y))$ can be expressed in
terms of partial derivatives of $P$.

Using the $C^3$-Lohner algorithm \cite{WZn} applied to the system
(\ref{eq:rosslerSystem}) with a range of parameter values $\bar A$
and an initial condition $Y\times \bar Z$ we computed bounds  of
partial derivatives of Poincar\'e map $P$ up to third order and an
estimation for $\frac{\partial^3 G}{\partial y^3}(\bar A\times Y)$
and $\frac{\partial^2 G}{\partial y^2}(\bar A\times Y)$. The
estimates (\ref{eq:d3gEstimation}) and (\ref{eq:d2gEstimation}) are
an interval enclosures of the estimates obtained in each of $30$
steps.

\end{itemize}

We used $6$-th order Taylor method with the time step $0.04$,
both, to verify the existence of $z(a,y)$ and to compute higher
order partial derivatives of $P$.
\qed

\subsection{The existence of period doubling bifurcation for $P$.}
\label{subsec:firstBif}
\begin{lemma}\label{lem:firstBifEstimations}
For $A=[a_1,a_2]$, $Y=[y_1,y_2]$ the following estimations hold true
\begin{eqnarray}
\frac{\partial G}{\partial y}(a_2,y_{fp}(a_2))&\in&
 [-1.2916325,-1.2916323]\cdot10^{-6} \label{eq:dgyAt_yfp}\\
G(a_2,y_1) &\in&[-1.15,-1.07]\cdot 10^{-13}\\
G(a_2,y_2) &\in&[5.2,5.21]\cdot 10^{-12}\label{eq:Gat_a2}\\
\frac{\partial^2 G}{\partial y\partial a}(A\times Y) &\subset&
[-2.421398492231531,-0.278863623843693] \label{eq:d2g_dydaEstimation}\\
\frac{\partial G}{\partial y}(\{a_1\}\times Y)&\subset& [0.83,16.87]
\cdot 10^{-6}
\end{eqnarray}

\end{lemma}
\textbf{Proof:} The estimations have been obtained using
$C^0-C^1-C^2$-Lohner algorithms applied to the systems
(\ref{eq:rosslerSystem}) and (\ref{eq:extendedRosslerSystem}). The
verification of conditions (\ref{eq:dgyAt_yfp}--\ref{eq:Gat_a2})
required computations in interval arithmetics based on $150$-bit
mantissa floating points.

The settings of $C^0$-$C^2$-Lohner methods for the above
computations are listed in Table~\ref{tab:settings}.
\begin{table}
\begin{center}
\begin{tabular}{|c|c|c|c|c|}
    \hline
    \mbox{ } & order & step & grid & remarks\\
    \hline
        $\displaystyle\frac{\partial G}{\partial y}(a_2,y_{fp}(a_2))$& $14$ & $0.05$ & --  & 150-bit mantissa\\
    \hline
        $\displaystyle G(a_2,y_1)$ & $14$ & $0.05$ &  --  & 150-bit mantissa\\
    \hline
       $ \displaystyle G(a_2,y_2)$ & $14$ & $0.05$ &  -- &  150-bit mantissa\\
    \hline
     $\displaystyle\frac{\partial^2 G}{\partial y\partial a}(A\times
    Y)$ & $6$ & $0.05$ & $5\times 30$ & integration of (\ref{eq:extendedRosslerSystem})\\
    \hline
    $\displaystyle\frac{\partial G}{\partial y}(\{a_1\}\times Y)$
    & $10$ & $0.05$ & $1\times 16000$ & nonequal parts\\
    \hline
\end{tabular}
\end{center}
\caption{Parameters of the $C^0-C^2$-Lohner
algorithms.\label{tab:settings}}
\end{table}
\qed

\textbf{Proof of Theorem~\ref{thm:rosslerFirstPD}:} The assertion
follows from Theorems~\ref{thm:per-doubling},~\ref{thm:dynInfo} and
numerical
Lemmas~\ref{lem:fixedPointCurve},~\ref{lem:derEigenvalue},~\ref{lem:LiapSchmidtExistence},~\ref{lem:firstBifEstimations}.

Indeed, assumptions of Theorem~\ref{thm:per-doubling} has been
verified in
\begin{itemize}
    \item \textbf{A1} -- Lemma~\ref{lem:LiapSchmidtExistence}
    \item \textbf{A2} -- from Lemma~\ref{lem:fixedPointCurve} there
        exists a fixed point curve $(y_{fp},z_{fp})\colon A\to
        (y_1,y_2)$ and from Lemma~\ref{lem:LiapSchmidtExistence} it has
        form as desired in \textbf{A2}
    \item \textbf{A3} -- $0\notin \frac{\partial^3 G}{\partial y^3}(A\times
    Y)$ because of (\ref{eq:d3gEstimation}).

    From (\ref{eq:d2gEstimation}), (\ref{eq:d2g_dydaEstimation}) and
    (\ref{eq:yfpPrimeBound}) it follows that $0\notin \frac{\partial^2
    G}{\partial a\partial y}(A\times Y) + \frac{\partial^2 G}{\partial
    y^2}y_{fp}'(A)\cdot[0,1]$.

    Finally, Lemma~\ref{lem:firstBifEstimations} guarantees that the
    remaining assumptions of Lemma~\ref{lem:Gsol} with
    $\varepsilon_1=\varepsilon_\nu=+1$ and $[\delta_1,\delta_2]=[y_1,y_2]$.

\end{itemize}
Finally, from Lemma~\ref{lem:derEigenvalue} we see that the
assumptions about the spectrum of $DP_a(A)$ and an eigenvector
$v(a)$ as desired in Theorem~\ref{thm:dynInfo} are satisfied. \qed

\subsection{The existence of second period doubling bifurcation.} \label{subsec:secondbif}

 In Section~\ref{subsec:firstBif} we
gave a computer assisted proof that for some parameter value $\bar
a_1\in[2.83244,2.832446]$ period doubling bifurcation occurs for
$P_{\bar a_1}$. In this section we use similar arguments in order to
prove that $P_{\bar a_2}^2$ has period doubling bifurcation for some
$\bar a_2\in[3.83735812,3.837358168411]$.

Since the arguments used to prove the existence of second period
doubling bifurcation are the same as in the first period doubling
bifurcation we omit the details and we present only the sets and the
necessary estimates.

Define
\begin{eqnarray*}
    A_2 &=& [a_3,a_4] = [3.83735812,3.837358168411]\\
    Y_2 &=& [y_3,y_4] = [-1.1,1.1]\cdot 10^{-6}\\
    Z_2 &=& \frac{1}{3}Y_2\\
    u_2 &=& (-4.5003284169596655673,0.043136987520848421584)\\
    M_2 &=& \begin{bmatrix}
            0.99999908059259889903 & 0.82277742767392003653\\
            0.0013560287448822982113 & -0.56836370794614177182
            \end{bmatrix}
\end{eqnarray*}

The point $u_2$ is an approximate period two point for parameter
value $a_4$, and the columns of matrix $M_2$ are normalized
eigenvectors of $DP^2_{a_4}$, where the first column corresponds to
eigenvalue close to $-1$.

On the Poincar\'e section $\Pi$ we will use a coordinates $(y,z)=
M^{-1}_2(u-u_2)$, where $u$ denotes a point in cartesian
coordinates. In this subsection we will use only these coordinates.

\begin{theorem}\label{thm:secBif}
The Poincar\'e map $P^2_a$ has a period doubling bifurcation at some
point $(\bar a_2,\bar y_2,\bar z_2)\in \inte(A_2\times Y_2\times
Z_2)$.
\end{theorem}
The proof is a consequence of the following  lemmas (proved with
computer assistance)

\begin{lemma}\label{lem:secFixedPointCurve}
There exist function $(y_{per},z_{per})\colon A_2\to Y_2\times Z_2$
smooth of class $C^\infty$ such that for $(a,y,z)\in A_2\times
Y_2\times Z_2$ holds
\begin{equation*}
P^2_a(y,z)=(y,z) \qquad \text{iff}\qquad
(y,z)=(y_{per}(a),z_{per}(a))
\end{equation*}
and
\begin{equation*}
y_{per}'(A_2)\subset
[-0.36435039423614490328,-0.36419313389173590956]
\end{equation*}
\end{lemma}

\begin{lemma}\label{lem:sec_derEigenvalue}
Let $\lambda_1,\lambda_2:A\to\mathbb R$ be eigenvalues of
$DP^2_a(y_{per}(a),z_{per}(a))$ defined by similar formulas as in
Lemma~\ref{lem:derEigenvalue}. Let $v(a)$ be the normalized
eigenvector corresponding to eigenvalue $\lambda_1(a)$. Then
\begin{eqnarray*}
    \lambda_1(a_3)  & \in     & [-0.99999992011934590863,-0.99999992005484927837]\\
    \lambda_1(a_4)  & \in     & [-1.00000000000149573159618,-1.00000000000149573159615]\\
    \lambda_2(A_2)  & \subset & [-7.7304566166653588839,7.7302177026359675856] \cdot 10^{-5}\\
    \lambda_1'(A_2) & \subset & [-1.6554066232416912996,-1.6460891324715511974] \\
    \pi_y v(A_2)        &\subset  &
    \pm[0.99984657385734598822,1.0001534381577519284].
\end{eqnarray*}
\end{lemma}

\begin{lemma}
For all $(a,y)\in A_2\times Y_2$ there exists unique $z=z(a,y)\in
Z_2$ such that
\begin{equation*}
P_{a,z}^4(y,z)=z\qquad \text{iff}\qquad z=z(a,y)
\end{equation*}
and the map $z\colon A_2\times Y_2\to Z_2$ is smooth of class
$C^\infty$. Moreover, the map $G\colon A_2\times Y_2\to \mathbb R$
defined by
\begin{equation*}
G(a,y) = y - P^4_{a,y}(y,z(a,y))
\end{equation*}
satisfies
\begin{eqnarray*}
\frac{\partial^3 G}{\partial y^3}(A\times Y) &\subset & [11.780861336872181511,22.544626008881969881] \\
\frac{\partial^2 G}{\partial y^2}(A\times Y) &\subset &
[-0.12474597648618415136,0.12474408945310766494]
\end{eqnarray*}
\end{lemma}

\begin{lemma}
The following estimations hold true
\begin{eqnarray*}
\frac{\partial G}{\partial y}(a_4,y_{per}(a_4))&\in&
[-2.992,-2.991]\cdot 10^{-12}
\\
G(a_4,y_3) &\in&[-5.13,-5.11]\cdot 10^{-19}\\
G(a_4,y_4) &\in&[5.21,5.22]\cdot 10^{-19}\\
\frac{\partial^2 G}{\partial y\partial a}(A_2\times Y_2) &\subset&
[-4.354355790892265432,-2.244876361570084633] \\
\frac{\partial G}{\partial y}(\{a_3\}\times Y_2)&\subset&
[0.99,30.98]\cdot 10^{-8}\\
\end{eqnarray*}
\end{lemma}

Parameter settings of computations involved in proofs of the above
lemmas are listed in Table~\ref{tab:secBif}.
\begin{table}
\begin{center}
\begin{tabular}{|c|c|c|c|c|}
    \hline
        \mbox{ } & order & step & grid & remarks\\
    \hline
        $\displaystyle\frac{\partial^3 G}{\partial y^3}(A_2\times
        Y_2)$ & $10$ & $0.04$ & $150\times 1$ & --\\
    \hline
        $\displaystyle\frac{\partial^2 G}{\partial y^2}(A_2\times
        Y_2)$ & $10$ & $0.04$ & $150\times 1$ & --\\
    \hline
        $\displaystyle\frac{\partial G}{\partial y}(a_4,y_{per}(a_4))$& $25$ & $0.05$ & --  & 150-bit mantissa\\
    \hline
        $\displaystyle G(a_4,y_3), G(a_4,y_4)$ & $25$ & $0.05$ &  --  & 150-bit mantissa\\
    \hline
     $\displaystyle\frac{\partial^2 G}{\partial y\partial a}(A_2\times
    Y_2)$ & $10$ & $0.05$ & $150\times 1$ & integration of (\ref{eq:extendedRosslerSystem})\\
    \hline
    $\displaystyle\frac{\partial G}{\partial y}(\{a_3\}\times Y_2)$
    & $14$ & $0.05$ & $1\times 10000$ & -- \\
    \hline
\end{tabular}
\end{center}
\caption{Parameters of the $C^0-C^3$-Lohner algorithms in the proof
of the existence of second period doubling
bifurcation.\label{tab:secBif}}
\end{table}

\section{Continuation of bifurcation diagram}
\label{sec:cont}

In the previous sections we proved that the map $P_a$ has period
doubling bifurcations for parameter values $\bar a_1\in [a_1,a_2]=A$
and $\bar a_2\in[a_3,a_4]=A_2$  in sets $Y \times Z$ and $Y_2 \times
Z_2$, respectively.

Our goal now is to connect these bifurcations with the curve of
period two points. More precisely, we prove the following result

\begin{theorem}\label{thm:connectingBranch}
There exists a continuous curve
\begin{equation*}
(y_{per},z_{per})\colon (\bar a_1,\bar a_2]\to \mathbb R^2
\end{equation*}
of period two points for $P_a$. Moreover,

\begin{eqnarray*}
(y_{per}(a),z_{per}(a)), P_a(y_{per}(a),z_{per}(a)) &\in& Y\times Z\quad \text{for } \bar
a_1<a\leq a_2\\
(y_{per}(a),z_{per}(a)), P_a^2(y_{per}(a),z_{per}(a)) &\in& Y_2\times Z_2\quad \text{for } a_3\leq
a\leq \bar a_2.
\end{eqnarray*}
Therefore curve $(y_{per},z_{per})$ connects the two bifurcation
points for $a=\bar{a}_1$ and $a=\bar{a}_2$.
\end{theorem}

The proof of the existence of a branch of period two points for $P$
consists of the following steps.
\begin{itemize}
\item[1.] the existence of continuous curve of period two points on
intervals $(\bar a_1,a_2]$ and $[a_3,\bar a_2]$ is a consequence of
Theorem~\ref{thm:rosslerFirstPD} and Theorem~\ref{thm:secBif},
respectively.
\item[2.] for parameter values slightly above $a_2$, $a_2 < a \leq \tilde{a}$, with $\tilde{a}-a_2$ small,  we extend
this curve using Lemma~\ref{lem:c3continuation}, which requires some
$C^3$ estimates (hence it is demanding computationally).
\item[3.] for parameters far from $a_2$ up to $a_3$, i.e. $\tilde{a} < a \leq a_3$, we verify the
existence of period two point curves using Krawczyk method
(Theorem~\ref{thm:KrawczykMethod}), which requires only $C^1$
estimates.
\item[4.] Since we use different methods for proving the
existence of segments of period two points curve  over some intervals in $[\bar{a}_1,\bar{a}_2]$ it is necessary
to verify that these segments can be glued to produce continuous curve.
\end{itemize}
At first, it appears that step 2, requiring costly $C^3$
computations, is not necessary, because in step 3 we can consider
also points close to $a_2$ using $C^1$ computations. But it turns
out that, while in principle possible, this approach may require
very large computation times, because the hyperbolicity is very
weak there, due  to the fact that one eigenvalue of $P^2_{a_2}$ is
very close to $1$.

To deal with this problem we used Lemma~\ref{lem:c3continuation}
to prove that for parameter values slightly above $a_2$ there
exists a continuous branch of period two points.
Algorithm~\ref{alg:d3g} is designed to verify assumptions of
Lemma~\ref{lem:c3continuation}. In Lemma~\ref{lem:corralg1} we
prove its correctness.

\begin{definition}
Let $U\subset \mathbb R^n$ be a bounded set. We say that $\mathcal
G\subset 2^{\mathbb R^n}$ is a grid of $U$ if
\begin{enumerate}
\item $\mathcal G$ is a finite set and each $G\in\mathcal G$ is a
closed set
\item $U\subset \bigcup_{G\in\mathcal G}G$
\end{enumerate}
\end{definition}

In our algorithms, which will be presented below, we always use
grids consisting of interval sets, i.e. sets which are cartesian
products of intervals, most of the time uniform grids, which are
defined as follows.

\begin{definition}
Let $Y=\Pi_{i=1}^n Y_i$, where $Y_i=[a_i,b_i]$ for $a_i \leq b_i$ and
let $(g_1,\dots,g_n) \in \mathbb{Z}_+$. We define \emph{a (uniform)  $g_1 \times g_2 \times \dots \times g_n$-grid for
  $Y$} denoted by $\mathcal{G}(g1_,\dots,g_n,Y)$ as follows.

For any $(j_1,\dots,j_n) \in \mathbb{Z}_+$, such that $j_i \leq
g_i $ we set
\begin{equation}
  g_{j_1,\dots,j_n}= \Pi_{i=1} \left[a_i, a_i + j_i \cdot \frac{b_i - a_i}{g_i}\right].
\end{equation}
Then $\mathcal{G}(g1_,\dots,g_n,Y)$ is a collection of all $g_{j_1,\dots,j_n}$.
\end{definition}

\begin{algorithm}\linesnumbered\label{alg:d3g}
 \caption{verification of assumptions of Lemma~\ref{lem:c3continuation}}
 \KwData{$[\nu_1,\nu_2]$ - an interval,
         $g_1, g_2, g_3, g_x$ - integers,
         $t\in(0,1)$ - float number,
         $X\times Y$ - a convex, compact set,
         $f_\nu$ - parameterized family of maps}
 \KwResult{If algorithms stops and does not throw an exception then assumptions of Lemma~\ref{lem:c3continuation} are satisfied}
 \Begin{
  $\mathcal G_1\longleftarrow$ $g_1\times g_x$-grid for $[\nu_1,\nu_2]\times X$\;
  \ForEach{$\bar\nu\times\bar X\in\mathcal G_1$}{
    $y(\bar\nu,\bar X) \longleftarrow \mathrm{IntervalNewtonOperator}(
        \Id_Y-\pi_Y\circ f^2_{\bar\nu}(\bar X,\cdot),Y,\mathrm{center}(Y)
        )$\;
    \If{\textbf{not} $y(\bar\nu,\bar X)\subset\inte Y$ }{
        {\textbf{throw} Liapunov-Schmidt reduction not verified}
    }
    \If{\textbf{not} $\frac{\partial^3 G}{\partial x^3}(\bar \nu,\bar X)>0$ }{
        {\textbf{throw} condition (\ref{eq:d3g}) is not satisfied}
    }
  }
  \If{(\textbf{not} $G(\nu_2,\min(X))<0$) \textbf{or} (\textbf{not} $G(\nu_2,\max(X))>0$)}{
    \textbf{throw} condition (\ref{eq:opositeSigns}) is not satisfied
    }
  $\mathcal G_2\longleftarrow$  $g_2$-grid for  $[\nu_1,\nu_2]$\;
  \ForEach{$\bar\nu\in\mathcal G_2$}{
    \If{(\textbf{not} $\frac{\partial G}{\partial\nu}(\bar\nu,x_1)>0$)
    \textbf{or}
    (\textbf{not} $\frac{\partial G}{\partial\nu}(\bar\nu,x_2)<0$)
    }
    {
    \textbf{throw} condition (\ref{eq:opositeSigns}) is not satisfied
    }
  }
  $\mathcal G_3\longleftarrow$  $g_3$-grid for  $[\nu_1,\nu_2]$\;
  \ForEach{$\bar\nu\in\mathcal G_3$}{
    $(\bar X,\bar Y) \longleftarrow \mathrm{IntervalNewtonOperator}(\Id-f_{\bar \nu},X\times Y,\mathrm{center}(X\times
    Y))$\;
    \If{\textbf{not} $(\bar X,\bar Y)\subset\inte(X\times Y)$}{
    \textbf{throw} fixed points curve not verified
    }
    $x_+\longleftarrow (1-t)\min(X) + t\min(\bar X)$\;
    $x_-\longleftarrow (1-t)\max(X) + t\max(\bar X)$\;
    \If{(\textbf{not} $G(\min(\bar\nu),x_+)>0$) \textbf{or} (\textbf{not} $G(\min(\bar\nu),x_-)<0$)}{
     \textbf{throw} condition (\ref{eq:xMinus}) or (\ref{eq:xPlus}) is not satisfied
    }
    \If{(\textbf{not} $\frac{\partial G}{\partial\nu}(\bar\nu,x_+)<0$)
    \textbf{or}
    (\textbf{not} $\frac{\partial G}{\partial\nu}(\bar\nu,x_-)>0$)
    }
    {
    \textbf{throw} condition (\ref{eq:xMinus}) or (\ref{eq:xPlus}) is not
    satisfied\;
    }
  }
 }
\end{algorithm}
\begin{lemma}
\label{lem:corralg1} If Algorithm~\ref{alg:d3g} is called with its
arguments $\nu_1,\nu_2,g_1,g_2,g_3,g_x,t$, $X\times Y$ and $f_\nu$
and it does not throw an exception then the assumptions of
Lemma~\ref{lem:c3continuation} are satisfied for $f_\nu$,
$\nu\in[\nu_1,\nu_2]$ on $X\times Y$.
\end{lemma}
\textbf{Proof:} The assumption about existence of fixed point curve
is verified in lines 15--19 since for all $\nu\in[\nu_1,\nu_2]$ the
Interval Newton Operator satisfies assumptions op
Theorem~\ref{thm:NewtonMethod}.

The existence of Liapunov-Schmidt reduction together with condition
(\ref{eq:d3g}) is verified in lines 2--8. In lines 4--6 we see that
$y(\nu,x)$ which solves equation $y-\pi_y(f^2_\nu(x,y))$ is unique
for fixed $\nu$, therefore by the implicit function theorem
$y(\nu,x)$ is smooth and we can compute map $G$ and its partial
derivatives.

In lines 9--10 we verify that $G(\nu_2,\min(X))<0$ and
$G(\nu_2,\max(X))>0$. Since $\frac{\partial G}{\partial
\nu}(\nu,\min(X))>0$ and $\frac{\partial G}{\partial
\nu}(\nu,\max(X))<0$ (lines 11--14) we see that for
$\nu\in[\nu_1,\nu_2]$ holds $G(\nu,\min(X))<0$ and
$G(\nu,\max(X))>0$. Therefore (\ref{eq:opositeSigns}) holds true.

Finally, in lines 15--25 we verify conditions
(\ref{eq:xMinus}--\ref{eq:xPlus}). Again we verify that for an
element of grid $\bar\nu$ it holds $G(\min(\bar\nu),x_+)>0$ and
$G(\min(\bar\nu),x_-)<0$. This together with $\frac{\partial
G}{\partial \nu}(\bar\nu,x_+)>0$ and $\frac{\partial G}{\partial
\nu}(\bar\nu,x_-)<0$ proves (\ref{eq:xMinus}--\ref{eq:xPlus}).

\qed

As was mentioned earlier  Algorithm~\ref{alg:d3g} is used to prove
the existence of period two points curve for parameter values
slightly above the first bifurcation - i.e. for parameters close
to $a_2$ $G$ has three solutions close to one another. For these
parameter values we found difficulties with verifying the
existence of period-two points curve using $C^1$-computations,
only.

For parameter values away from the bifurcation, where all
eigenvalues of periodic orbit are well separated from the unit
circle, we use  Algorithm~\ref{alg:Krawczyk} based on the Newton
interval method and the Krawczyk method, both requiring only
$C^1$-computations,   and which verifies the existence of only one
branch of period two points for $P_a$.

Before we present this algorithm we need to introduce some
notations. Let $\bar\Pi =\{(x,y,z)\in\mathbb R^3 : x=0\}$ be a
Poincar\'e section for (\ref{eq:rosslerSystem}) and
$\overline{P}_a\colon \bar\Pi\oarrow\bar\Pi$ be corresponding
Poincar\'e map for a system with parameter value $a$. Notice that
the trajectory can intersect $\bar\Pi$ at a point
$(y,z)\in\bar\Pi$ for which $x'=-y-z$ is positive or negative (if
it is equal to zero the Poincar\'e map is not defined). Hence we
have $\bar P^2_a|_\Pi = P_a$, where $P_a$ is the Poincar\'e map
for section $\Pi=\{(x,y,z)\in\mathbb R^3\ |\ x=0,x'>0\}$ and
therefore period two points for $P_a$ correspond to period four
points for $\bar P_a$. Let us define a map
$F_a:\bar\Pi^4\to\mathbb R^8$ by
\begin{equation*}
F_a\begin{bmatrix}
    (y_1,z_1)\\
    (y_2,z_2)\\
    (y_3,z_3)\\
    (y_4,z_4)
\end{bmatrix} =
\begin{bmatrix}
    (y_2,z_2) - \bar P_a(y_1,z_1)\\
    (y_3,z_3) - \bar P_a(y_2,z_2)\\
    (y_4,z_4) - \bar P_a(y_3,z_3)\\
    (y_1,z_1) - \bar P_a(y_4,z_4)\\
\end{bmatrix}
\end{equation*}

Algorithm~\ref{alg:Krawczyk} was used to verify the existence of a
continuous branch of period two points for $P_a$ for $a$ belonging
to some interval. Sets $X_i \times Y_i$ give the size of the
neighborhood around a candidate periodic orbit on section
$\overline{\Pi}$. Lines 4 to 8 constitute a heuristic part and their
task is to find a good  candidate.

\begin{algorithm}\linesnumbered\label{alg:Krawczyk}
\caption{verification the existence of period two points branch.}
 \KwData{$[a_*,a^*]$ - an interval, $Y_i\times Z_i$, $i=1,2,3,4$ - convex, compact sets, $g$ - an integer}
 \KwResult{if algorithm stops and does not throw an exception then there exists a continuous branch of period two points for $P_a$ for parameter values $a\in[a_*,a^*]$}
 \Begin{
    $\mathcal G \longleftarrow$ $g$-grid for $[a_*,a^*]$\;
    \ForEach{$\bar a\in\mathcal G$}{
        $a\longleftarrow\mathrm{center}(\bar a)$\;
        $u_1=(y_1,z_1)\longleftarrow$ find approximate period two point for $P^2_a$ using
        standard Newton method\;
        $u_2=(y_2,z_2)\longleftarrow \bar P_a(y_1,z_1)$\;
        $u_3=(y_3,z_3)\longleftarrow \bar P_a(y_2,z_2)$\;
        $u_4=(y_4,z_4)\longleftarrow \bar P_a(y_3,z_3)$\;
        $C\longleftarrow$ compute approximate value of
        $DF_a(u_1,u_2,u_3,u_4)$\;
        \If{$C$ is singular}{$C\longleftarrow \Id$\;}
        $U\longleftarrow(u_1+Y_1\times Z_1,u_2+Y_2\times Z_2,u_3+Y_3\times Z_3,u_4+Y_4\times
        Z_4)$\;
        $u\longleftarrow(u_1,u_2,u_3,u_4)$\;
        $K_{\bar a}=(k_{1,\bar a},k_{2,\bar a},k_{3,\bar a},k_{4,\bar a})\longleftarrow \mathrm{IntervalKrawczykOperator}(F_{\bar a},C^{-1},U,u)$\;
        \If{\textbf{not} $K_{\bar a}\subset\inte U$}{
            \textbf{throw} cannot verify the existence of period two
            point\;
        }
        \If{$k_{1,\bar a}\cap k_{3,\bar a}\neq\emptyset$}{
            \textbf{throw} the unique fixed point for $P^2_{\bar a}$ in $k_{1,\bar a}$ is not necessary
            period two point for $P_{\bar a}$\;
        }
    }
    $B\longleftarrow\bigcup_{\bar a\in\mathcal G} \bar a\times k_{1,\bar
    a}$\;
    \If{$B$ is not connected}{
        \textbf{throw} cannot verify if branch of fixed point
        curve is continuous on interval $[a_*,a^*]$\;
    }
 }
\end{algorithm}

\begin{lemma}
\label{lem:algKrawczykOk}
If Algorithm~\ref{alg:Krawczyk} is called with its arguments
$[a_*,a^*]$, $Y_i\times Z_i$, $i=1,2,3,4$ and $g$ and does not throw
an exception then there exists a continuous curve
$(y_{per},z_{per})\colon [a_*,a^*]\to \Pi$ such that
$(y_{per}(a),z_{per}(a))$ is period two point for $P_a$.
\end{lemma}
\textbf{Proof:} The existence of fixed point for $P^2_a$ for all
$a\in[a_*,a^*]$ is verified in lines 12--16. Lines 17--18 guarantee
that this is a period two point for $P_a$, in fact a unique one in
$U$.

Uniqueness implies continuity on each $\bar{a} \in \mathcal{G}$
and due to uniqueness and connectedness of the set $B$ defined in
line 19 we see that they agree on boundaries of $\bar{a}$.

\qed

\textbf{Proof of Theorem~\ref{thm:connectingBranch}:}

 The existence of continuous curve of period two points on
intervals $(\bar a_1,a_2]$ and $[a_3,\bar a_2]$ is a consequence of
Theorem~\ref{thm:rosslerFirstPD} and Theorem~\ref{thm:secBif},
respectively. Let $a_s=2.8329$.  For parameter values $[a_2,a_s]$ we verify the
existence of period two points branch using Algorithm~\ref{alg:d3g}
and for parameter values $a\in[a_s,a_3]$ we use
Algorithm~\ref{alg:Krawczyk}.

We have ran Algorithm~\ref{alg:d3g} five times with parameters
listed in Table~\ref{tab:alg3dData} (in each case the map is
$P_a$). Since in each case the algorithm  had stopped and did not
throw an exception we conclude that in each interval of parameters
listed in Table~\ref{tab:alg3dData} there exist two continuous
curves $c_1(a)$, and $c_2(a)$ of period two points. Sets
$X_i\times Y_i$ listed in Table~\ref{tab:alg3dData} are chosen so
that
\begin{equation}\label{eq:monSequence}
Y\times Z\subset X_1\times Y_1\subset \cdots \subset X_5\times Y_5,
\end{equation}
where $Y\times Z$ is the set used in the proof of the existence of
first period doubling bifurcation.
 Observe also, that since we know
that for $a_2$ holds $P_{a_2}(c_1(a_2))=c_2(a_2)$ and
$c_1(a_2),c_2(a_2)\in Y\times Z\subset X_1\times Y_1$ from
Lemma~\ref{lem:c3continuation} we obtain that
$\left\{c_1(a),c_2(a)\right\}$ is period two orbit for $P_a$, for
$a\in[a_2,a_s]$, i.e. the whole interval of parameters covered by
intervals listed in the first columns in
Table~\ref{tab:alg3dData}. The uniqueness of period two orbit
together with condition (\ref{eq:monSequence}) implies that the
curves are continuous on $(\bar a_1,a_s]$.

One can see that the total number of initial values for which we
need compute third order derivatives of $G$,which  is equal to the
sum of $g_1g_x$ over all rows in Table~\ref{tab:alg3dData}, is equal
to $19350$. The total time of computation of this step is ten hours
on the Pentium IV, 3GHz processor.

\begin{table}
\begin{center}
\begin{tabular}{|c|c|c|c|c|c|c|}
    \hline
    $[\nu_1,\nu_2]$ & $g_1$ & $g_2$ & $g_3$ & $g_x$ & $t$ & $(X_i\times
    Y_i)\cdot 10^{-4}$\\
    \hline
    $[a_2,2.8325]$ & $6$ & $70$ & $600$ & $150$ & $0.88$ &
    $[-99.0216,97.95]\times[-4,4]$\\
    \hline
    $[2.8325,2.8326]$ & $10$ & $80$ & $300$ & $270$ & $0.7$ &
    $[-166.9,163.9]\times[-4,4]$\\
    \hline
    $[2.8326,2.8327]$ & $12$ & $60$ & $150$ & $320$ & $0.6$ &
    $[-214.59,209.62]\times[-4,4]$\\
    \hline
    $[2.8327,2.8328]$ & $17$ & $50$ & $120$ & $330$ & $0.5$ &
    $[-253.8,246.8]\times[-8,8]$\\
    \hline
    $[2.8328,2.8329]$ & $18$ & $50$ & $100$ & $350$ & $0.5$ &
    $[-287.757,279]\times[-8,8]$\\
    \hline
\end{tabular}
\end{center}
\caption{Parameters of the
Algorithm~\ref{alg:d3g}.\label{tab:alg3dData}}
\end{table}
We have run Algorithm~\ref{alg:Krawczyk} with $74$
different arguments listed in Table~\ref{tab:algKrawczykData}.
 We have chosen the parameters of the
Algorithm~\ref{alg:Krawczyk} so that such $74$ intervals
$[(a_*)_i,(a^*)_i]$, $i=1,\ldots,74$ cover the interval
$[a_s,a_3]$. Notice also, that for parameters $a$ closer to $a_s$
we need larger values of $g$ since the hyperbolicity close to
$a_2$ is very weak. The total number of subintervals used  to
cover  interval $[a_s,a_3]$ is $614450$. In fact this is the
longest part of the numerical proof. The total time of computation
of this step is 53 hours on the Pentium IV, 3GHz processor.  Since
in each case Algorithm~\ref{alg:Krawczyk} stops and does not throw
an exception we conclude that on each subinterval
$[(a_*)_i,(a^*)_i]$ there exists continuous branch of period two
points. We need to show that these curves glue continuously at
$a^*_i$'s.  In fact, this algorithm returns an upper bound for
this period two points branch which is of the form
\begin{equation}
    B=\bigcup_{i=1}^{74}B_i \label{eq:defW-alg2}
\end{equation}
where $B_i$'s are defined in line 19 of
Algorithm~\ref{alg:Krawczyk}. We verified that $B$ is connected -
this together with an information that for fixed $a\in[a_s,a_3]$
there exists a unique period two point $(y_{per}(a),z_{per}(a))$
such that $(a,y_{per}(a),z_{per}(a))\in B$ implies that the curve
$(y_{per}(a),z_{per}(a))$ is continuous on $[a_s,a_3]$.

There remains to show the continuity of fixed point branch for
parameter values $a=a_s$ and $a=a_3$. For $a=a_s$ we know that
there exist period two points $c_1(a_s), c_2(a_s)$ which belongs
to the last set listed in Table~\ref{tab:alg3dData}, i.e.
\begin{equation*}
c_1(a_s),c_2(a_s)\in W_1=u_0 +
M\cdot([-287.757,279]\times[-8,8])\cdot 10^{-4}
\end{equation*}
where $u_0$ and $M$ define coordinate system close to first period
doubling bifurcation and are defined in
(\ref{eq:firstBifPoint}--\ref{eq:firsBifCoordSystem}). On the other
hand the estimation for period two point resulting from Krawczyk
method used in Algorithm~\ref{alg:Krawczyk} is
\begin{eqnarray*}
W_2 &=& (W_2^1,W_2^2)\\
W_2^1 &=&[-4.7668051788293892557,-4.7667832743925968586]\\
W_2^2 &=&[0.052543190547910088861,0.052543238016254205369]
\end{eqnarray*}
One can verify that $W_2\subset W_1$ which obviously means that a
period two point $(y_{per}(a_s),z_{per}(a_s))\in W_2$ resulting from
the Krawczyk method and Algorithm~\ref{alg:Krawczyk} is one of the
points $c_1(a_s),c_2(a_s)$ resulting from Algorithm~\ref{alg:d3g}.
Hence, the curve of period two points is continuous at $a=a_s$.

In a similar way we verified continuity at $a=a_3$. From the
Krawczyk method used in Algorithm~\ref{alg:Krawczyk} we know that
$(y_{per}(a_3),z_{per}(a_3))$ is a unique period two point in the
set
\begin{eqnarray*}
W_3 &=& (W_3^1,W_3^2)\\
W_3^1 &=&[-4.5010116820607413146,-4.4996232549240025023]\\
W_3^2 &=&[0.043134233933640332703,0.043140290681655812932]
\end{eqnarray*}
On the other hand from Theorem~\ref{thm:secBif} we know that for
$a=a_3$ period two point belongs to the set
\begin{equation*}
W_4 = u_2 + M_2\cdot(Y_2\times Z_2)
\end{equation*}
where $M_2$, $u_2$, $Y_2,Z_2$ define the set on which we verify the
existence of second period doubling bifurcation. One can verify that
$W_4\subset W_3$ which proves that the branch of period two points
is continuous at $a=a_3$.

\qed

\begin{center}
\tablefirsthead{%
    \hline
    \multicolumn{4}{|c|}{$U=[-1,1]\cdot(1100,3,1000,3000,1100,3,1000,3000)\cdot
10^{-8}$}\\
    \hline
  $i$ & $[a_*,a^*]$ & $g$ & $Y_1\times Z_1\times Y_2\times Z_2\times Y_3\times Z_3\times Y_4\times Z_4$\\
  \hline
}%
\tablehead{%
\hline
\multicolumn{4}{|l|}{\small\sl Table~\ref{tab:algKrawczykData} continued from previous page}\\
\hline
    \hline
  $i$ & $[a_*,a^*]$ & $g$ & $Y_1\times Z_1\times Y_2\times Z_2\times Y_3\times Z_3\times Y_4\times Z_4$\\
  \hline
}%
\tabletail{%
\hline
\multicolumn{4}{|l|}{\small\sl Table~\ref{tab:algKrawczykData} continued on the next page}\\
\hline
}%
\tablelasttail{
    \hline
}%
\bottomcaption{Parameters of Algorithm~\ref{alg:Krawczyk}. The
initial set $U$ is defined in the first line of the table.
 \label{tab:algKrawczykData}}

\begin{supertabular}{|c|c|c|c|}
 \hline $1$ & $[2.8329,2.83291]$ & $14000$ & $U$\\
 \hline$2$ & $[2.83291,2.83292]$ & $13500$ & $U$\\
 \hline$3$ & $[2.83292,2.83293]$ & $12200$ & $3U$\\
 \hline$4$ & $[2.83293,2.83294]$ & $11250$ & $3U$\\
 \hline$5$ & $[2.83294,2.83295]$ & $10400$ & $3U$\\
 \hline$6$ & $[2.83295,2.83296]$ & $9650$ & $3U$\\
 \hline$7$ & $[2.83296,2.83297]$ & $9050$ & $3U$\\
 \hline$8$ & $[2.83297,2.83298]$ & $8500$ & $3U$\\
 \hline$9$ & $[2.83298,2.83299]$ & $8000$ & $3U$\\
 \hline$10$ & $[2.83299,2.833]$ & $7550$ & $3U$\\
 \hline$11$ & $[2.833,2.83301]$ & $7200$ & $3U$\\
 \hline$12$ & $[2.83301,2.83302]$ & $6800$ & $3U$\\
 \hline$13$ & $[2.83302,2.83303]$ & $6500$ & $3U$\\
 \hline$14$ & $[2.83303,2.83304]$ & $6200$ & $3U$\\
 \hline$15$ & $[2.83304,2.83305]$ & $6000$ & $3U$\\
 \hline$16$ & $[2.83305,2.83306]$ & $5700$ & $3U$\\
 \hline$17$ & $[2.83306,2.83307]$ & $5500$ & $3U$\\
 \hline$18$ & $[2.83307,2.83308]$ & $5300$ & $3U$\\
 \hline$19$ & $[2.83308,2.83309]$ & $5100$ & $3U$\\
 \hline$20$ & $[2.83309,2.8331]$ & $4900$ & $3U$\\
 \hline$21$ & $[2.8331,2.83311]$ & $4800$ & $3U$\\
 \hline$22$ & $[2.83311,2.83312]$ & $4600$ & $3U$\\
 \hline$23$ & $[2.83312,2.83313]$ & $4450$ & $3U$\\
 \hline$24$ & $[2.83313,2.83314]$ & $4300$ & $3U$\\
 \hline$25$ & $[2.83314,2.83315]$ & $4150$ & $3U$\\
 \hline$26$ & $[2.83315,2.83316]$ & $4050$ & $3U$\\
 \hline$27$ & $[2.83316,2.83317]$ & $3950$ & $3U$\\
 \hline$28$ & $[2.83317,2.83318]$ & $3850$ & $3U$\\
 \hline$29$ & $[2.83318,2.83319]$ & $3750$ & $3U$\\
 \hline$30$ & $[2.83319,2.8332]$ & $3650$ & $3U$\\
 \hline$31$ & $[2.8332,2.8333]$ & $36000$ & $3U$\\
 \hline$32$ & $[2.8333,2.8334]$ & $29000$ & $3U$\\
 \hline$33$ & $[2.8334,2.8335]$ & $24000$ & $3U$\\
 \hline$34$ & $[2.8335,2.8336]$ & $20000$ & $3U$\\
 \hline$35$ & $[2.8336,2.8337]$ & $18000$ & $3U$\\
 \hline$36$ & $[2.8337,2.8338]$ & $16000$ & $3U$\\
 \hline$37$ & $[2.8338,2.8339]$ & $14000$ & $3U$\\
 \hline$38$ & $[2.8339,2.834]$ & $13000$ & $3U$\\
 \hline$39$ & $[2.834,2.8345]$ & $59000$ & $3U$\\
 \hline$40$ & $[2.8345,2.835]$ & $42000$ & $3U$\\
 \hline$41$ & $[2.835,2.8355]$ & $17500$ & $15U$\\
 \hline$42$ & $[2.8355,2.836]$ & $11000$ & $15U$\\
 \hline$43$ & $[2.836,2.8365]$ & $8000$ & $15U$\\
 \hline$44$ & $[2.8365,2.837]$ & $6200$ & $15U$\\
 \hline$45$ & $[2.837,2.8372]$ & $2800$ & $30U$\\
 \hline$46$ & $[2.8372,2.8375]$ & $3600$ & $30U$\\
 \hline$47$ & $[2.8375,2.838]$ & $5000$ & $30U$\\
 \hline$48$ & $[2.838,2.8385]$ & $3600$ & $30U$\\
 \hline$49$ & $[2.8385,2.839]$ & $2900$ & $30U$\\
 \hline$50$ & $[2.839,2.8395]$ & $2400$ & $30U$\\
 \hline$51$ & $[2.8395,2.84]$ & $2100$ & $30U$\\
 \hline$52$ & $[2.84,2.841]$ & $3700$ & $30U$\\
 \hline$53$ & $[2.841,2.842]$ & $3000$ & $30U$\\
 \hline$54$ & $[2.842,2.843]$ & $2500$ & $30U$\\
 \hline$55$ & $[2.843,2.844]$ & $2200$ & $30U$\\
 \hline$56$ & $[2.844,2.845]$ & $2000$ & $30U$\\
 \hline$57$ & $[2.845,2.846]$ & $1700$ & $30U$\\
 \hline$58$ & $[2.846,2.848]$ & $3100$ & $30U$\\
 \hline$59$ & $[2.848,2.85]$ & $2600$ & $30U$\\
 \hline$60$ & $[2.85,2.86]$ & $11000$ & $30U$\\
 \hline$61$ & $[2.86,2.87]$ & $6500$ & $30U$\\
 \hline$62$ & $[2.87,2.88]$ & $4700$ & $30U$\\
 \hline$63$ & $[2.88,2.89]$ & $3700$ & $30U$\\
 \hline$64$ & $[2.89,2.9]$ & $3100$ & $30U$\\
 \hline$65$ & $[2.9,2.95]$ & $4400$ & $150U$\\
 \hline$66$ & $[2.95,3]$ & $2500$ & $150U$\\
 \hline$67$ & $[3,3.1]$ & $3500$ & $150U$\\
 \hline$68$ & $[3.1,3.2]$ & $2500$ & $150U$\\
 \hline$69$ & $[3.2,3.3]$ & $2100$ & $150U$\\
 \hline$70$ & $[3.3,3.4]$ & $1800$ & $150U$\\
 \hline$71$ & $[3.4,3.5]$ & $1600$ & $150U$\\
 \hline$72$ & $[3.5,3.6]$ & $1400$ & $150U$\\
 \hline$73$ & $[3.6,3.7]$ & $1500$ & $150U$\\
 \hline$74$ & $[3.7,a_3]$ & $2400$ & $150U$\\
\hline
\end{supertabular}
\end{center}

\subsection{Technical data.}
In order to compute Poincar\'e maps $P$ and $P^2$ with their partial
derivatives we used the interval arithmetic \cite{IE,Mo}, the set
algebra and the $C^r$-Lohner algorithm \cite{WZn} developed at the
Jagiellonian University by the CAPD group \cite{CAPD}. The C++
source files of the program with an instruction how it should be
compiled and run are available at \cite{WI}.

All computations were performed with the Pentium IV, 3GHz
processor and 512MB RAM under Kubuntu Feisty Fawn linux with
gcc-4.1.1 and MS Windows XP Professional with gcc-3.4.4. The
computations took approximately three days. The main
time-consuming part (over 63 hours) is the verification of the
existence of connecting branch of period two points between first
and second bifurcation.

\end{document}